%% file: qcls.tex
\documentstyle[amssymb]{amsart}

\numberwithin{equation}{section}

\theoremstyle{plain}
\newtheorem{thm}{Theorem}[section]
\newtheorem{cor}[thm]{Corollary}
\newtheorem{prop}[thm]{Proposition}
\newtheorem{lem}[thm]{Lemma}

\theoremstyle{definition}
\newtheorem{defn}[thm]{Definition}

\def\ti#1{\mbox{\tiny $#1$}}
\def\scri#1{\mbox{\scriptsize $#1$}}
\def\ft#1{\mbox{\footnotesize $#1$}}

\def\bol#1{\boldsymbol{#1}}

\newcommand{\al}{\alpha}
\newcommand{\be}{\beta}
\newcommand{\ga}{\gamma}
\newcommand{\de}{\delta}
\newcommand{\vep}{\varepsilon}
\newcommand{\lam}{\lambda}
\newcommand{\Lam}{\Lambda}
\newcommand{\ol}{\overline}

\newcommand{\wt}{\widetilde}

\newcommand{\wtil}{\widetilde}
\newcommand{\st}{\stackrel}
\newcommand{\cl}{\cal} %
\newcommand{\fr}{\frak} %
\newcommand{\bd}{\bold} %

\newcommand{\SI}{\st{\circ}{\bd I}}

\newcommand{\tiJ}{\ti{\bd J}}
\newcommand{\tiK}{\ti{\bd K}}
\newcommand{\tiL}{\ti{\bd L}}
\newcommand{\tiS}{\ti{\bd S}}

\newcommand{\scJ}{\scri{\bd J}}
\newcommand{\scK}{\scri{\bd K}}
\newcommand{\scL}{\scri{\bd L}}
\newcommand{\sPi}{\st{\circ}{\Pi}}
\newcommand{\sDe}{\st{\circ}{\De}}
\newcommand{\sW}{\st{\circ}{W}}
\newcommand{\De}{\Delta}
\newcommand{\pr}{\prec}
\newcommand{\preq}{\preceq}

\newcommand{\re}{\De^{re}}
\newcommand{\pre}{\De^{re}_+}
\newcommand{\nre}{\De^{re}_-}
\newcommand{\im}{\De^{im}}
\newcommand{\pim}{\De^{im}_+}

\newcommand{\Lra}{\Longrightarrow}
\newcommand{\lra}{\longrightarrow}
\newcommand{\la}{\langle}
\newcommand{\ra}{\rangle}

\newcommand{\vPhi}{\varPhi}
\newcommand{\nab}{\nabla}

\newcommand{\br}{\mbox{\boldmath $r$}}
\newcommand{\bs}{\mbox{\boldmath $s$}}

\newcommand{\bul}{\bullet}

\title[\null]{The classification of convex orders \\
 on affine root systems}

\begin{document}
\maketitle

\begin{center}
 Ken Ito \\
 Graduate School of Mathematics, Nagoya University, \\
 Chikusa-ku, Nagoya, 464-8602, Japan \\
 e-mail: ken-ito{\char'100}math.nagoya-u.ac.jp
\end{center}

\vspace{0.5em}
\begin{quotation}
{\footnotesize
\begin{center}
 {\bf Abstract}
\end{center}

 We classify all total orders having a certain convex property
 on the positive root system of an arbitrary
 untwisted affine Lie algebra ${\frak g}$.
 Such total orders are called convex orders and are used to
 construct convex bases of Poincar\'e-Birkhoff-Witt type of
 the upper triangular subalgebra $U_q^+$ of
 the quantized enveloping algebra $U_q({\frak g})$.}
\end{quotation}

\vspace{0.5em}
\input{sec1}
\input{sec2}
\input{sec3}
\input{sec4}
\input{sec5}
\input{sec6}
\input{sec7}

\vspace{1em}
\begin{center}
 {\sc Acknowledgements}
\end{center}
 I thank Professor Akihiro Tsuchiya and Professor Takahiro Hayashi
 for constant help and precious advice.

\end{document}

%% file: sec1.tex
\section{Introduction}
 Let $\De$ be the root system of a Kac-Moody Lie algebra ${\fr g}$,
 $\De_+$ the set of all positive (resp. negative) roots relative to
 a root basis $\Pi=\{\,\al_i\,|\,i\in{\bd I}\,\}$ of $\De$,
 and $\pre$ (resp. $\pim$) the set of
 all positive real (resp. imaginary) roots.
 Denote by $s_i$ the simple reflection associated with
 $\al_i$ for each $i\in{\bd I}$, and by
 $W=\la{\,s_i\,|\,i\in{\bd I}\,}\ra$ the Weyl group of ${\fr g}$.
 Then $(W,S)$ is a Coxeter system with
 $S=\{\,s_i\,|\,i\in{\bd I}\,\}$ (cf. \cite{vK}).
 A total order $\preq$ on $\De_+$ is called a {\it convex order\/}
 if it satisfies the following condition:
\[ (\be,\ga)\in(\De_+)^2\setminus(\pim)^2, \;\; %
 \be\,\pr\,\ga,\;\; \be+\ga\in\De_+ \;\;\Lra\;\; %
 \be\,\pr\,\be+\ga\,\pr\ga. \]
 We call an infinite sequence
 $\bs=(\bs(p))_{p\in{\bd N}}\in S^{\bd N}$
 an {\it infinite reduced word\/} of $(W,S)$
 if the length of the element
 $z_{\bs}(p):=\bs(1)\cdots\bs(p)\in W$
 is $p$ for each $p\in{\bd N}$.

 The purpose of this paper is to classify
 all convex orders on $\De_+$ in case ${\fr g}$
 is an arbitrary untwisted affine Lie algebra.
 The motive of this study is related to
 the construction of {\it convex bases\/} of
 the upper triangular subalgebra $U_q^+$ of
 the quantized enveloping algebra $U_q({\fr g})$.
 Convex bases are Poincar\'e-Birkhoff-Witt type bases
 having a convex property concerning the ``$q$-commutator''
 of two ``$q$-root vectors'' of $U_q^+$.
 The convex property is useful to calculate values of
 the standard bilinear pairing between $U_q^+$ and
 the lower triangular subalgebra $U_q^-$, and is applied for
 explicit calculation of the universal $R$-matrix of
 $U_q({\fr g})$ (cf. \cite{KR}, \cite{KT}).
 By the way, each convex basis of $U_q^+$ is formed by
 monomials in certain $q$-root vectors of $U_q^+$
 multiplied in a predetermined convex order on $\De_+$.
 Thus convex orders on $\De_+$ are indispensable for
 the construction of convex bases of $U_q^+$.

 In case ${\fr g}$ is an arbitrary finite-dimensional simple
 Lie algebra, it is known that there exists a natural bijective
 mapping between the set of all convex orders on $\De_+$ and
 the set of all reduced expressions of the longest element of
 $W$ (cf. \cite{pP1}). In \cite{gL}, G. Lusztig constructed
 convex bases of $U_q^+$ associated with all reduced expressions
 of the longest element of $W$ by using a braid group action on
 $U_q({\fr g})$, so all convex bases of $U_q^+$
 have been constructed in this case.

 In case ${\fr g}$ is an arbitrary untwisted affine Lie algebra,
 convex orders on $\De_+$ are closely related to
 infinite reduced words of $(W,S)$.
 More precisely, each infinite reduced word naturally
 corresponds to a ``{\it 1-row type\/}'' convex order on
 an infinite proper subset of $\pre$ (see Section 3).
 Such infinite proper subsets of $\pre$ have
 interesting two convex properties, so we call them
 {\it infinite real biconvex sets\/}.
 Moreover, each convex order on $\De_+$ is made from each couple of
 ``{\it maximal\/}'' (infinite) real biconvex sets with
 convex orders which divides $\pre$ into two parts (see Theorem 7.3).
 In \cite{pP2}, P. Papi constructed all 1-row type convex orders on
 each maximal real biconvex set.
 In \cite{jB}, J. Beck constructed convex bases of $U_q^+$
 associated with convex orders on $\De_+$ which are made from a
 certain couple of maximal real biconvex sets with
 1-row type convex orders which divides $\pre$ into two parts.
 However, we expect that it is possible to generalize the Beck's
 construction, since we find that there exist several types of
 convex order called ``{\it n-row types}'' on each maximal
 real biconvex set on top of 1-row type (see Definition 3.1).
 In \cite{kI}, to analyze convex orders on
 maximal real biconvex sets, we studied detailed relationships
 between the set of all infinite reduced words and
 the set of all infinite real biconvex sets.

 In this paper, we show that the classification of
 convex orders on $\De_+$ is reduced to that of
 convex orders on maximal real biconvex sets (see Theorem 7.3),
 and classify and construct all convex orders on
 each maximal real biconvex set (see Theorem 7.9), and give
 a general method of the construction of convex orders on $\De_+$
 by summarizing all results (see Corollary 7.10).
 More precisely, we show that $n$-row type convex orders on
 each maximal real biconvex set $M$ naturally correspond
 to chains $C_1\subset\cdots\subset C_n$ consisting of
 infinite real biconvex sets $C_1,\dots,C_n$ such that
 $C_n=M$ and $\sharp(C_{i'}\setminus C_i)=\infty$ for $i<i'$,
 and we give a parameterization of the set of all such chains of
 infinite real biconvex sets (see Theorem 6.9).

 The text is organized as follows.
 In Section 2, we give the definition of biconvex sets
 for root systems with Coxeter group actions,
 and state several fundamental results.
 We next give the definition of infinite reduced words
 for each Coxeter system $(W,S)$, and introduce a natural
 correspondence between the set of all infinite reduced words of
 $(W,S)$ and the set of all infinite real biconvex sets.
 In Section 3, we prepare notation for total orders on sets,
 and give the definition of convex orders on
 $B\subset\De_+$ for each root system $\De$ of
 an arbitrary Coxeter system $(W,S)$,
 and introduce several relationships between convex orders,
 biconvex sets, and infinite reduced words.
 From Section 4 until Section 7, we treat only the case that
 ${\fr g}$ is an arbitrary untwisted affine Lie algebra.
 In Section 4, we prepare notation for untwisted affine cases
 and define a family of Coxeter systems $(W_{\tiJ},S_{\tiJ})$
 and its root systems $\De_{\tiJ}$.
 In Section 5, we describe detailed relationships between
 the set ${\cl W}_{\tiJ}^\infty$ of all infinite reduced words of
 $(W_{\tiJ},S_{\tiJ})$ and the set ${\fr B}_{\tiJ}^\infty$ of
 all infinite real biconvex sets in $\De_{\tiJ}$.
 In Section 6, we consider chains $C_1\subset\cdots\subset C_n$
 consisting of elements of ${\fr B}_{\tiJ}^\infty$.
 In Section 7, we classify all convex orders on $\De_+$ and give
 a general method of the construction of convex orders on $\De_+$.

%% file: sec2.tex
\section{Definitions of biconvex sets and reduced words}
 Let ${\bd R}$, ${\bd Z}$, and ${\bd N}$
 be the set of the real numbers, the integers,
 and the positive integers, respectively.
 Set ${\bd R}_{\geq a}:=\{b\in{\bd R}\,|\,b\geq a\}$,
 ${\bd R}_{>a}:=\{b\in{\bd R}\,|\,b>a\}$,
 and ${\bd Z}_{\geq a}:={\bd Z}\cap{\bd R}_{\geq a}$
 for each $a\in{\bd R}$.
 Denote by $\sharp U$ the cardinality of a set $U$.
 When $A$ and $B$ are subsets of $U$, write
 $A\,\dot{\subset}\,B$ or $B\,\dot{\supset}\,A$
 if $\sharp(A\setminus B)<\infty$, and write $A\,\dot{=}\,B$
 if $A\,\dot{\subset}\,B$ and $A\,\dot{\supset}\,B$.
 Then $\dot{=}$ is an equivalence relation on the power set of $U$.

\vspace{1em}
 Let $W$ be a group generated by a set $S$ of
 involutive generators (i.e. $s\neq1,\,s^2=1\;\forall\,s\in S$),
 and $(V,\De,\Pi)$ a triplet satisfying
 the following conditions R(i)--R(iv).

 R(i)\; It consists of a representation space $V$ of $W$ over
 ${\bd R}$, a $W$-invariant subset $\De\subset V\setminus\{0\}$
 which is symmetric (i.e. $\De=-\De$), and
 a subset $\Pi=\{\,\al_s\,|\,s\in S\,\}\subset\De$.

 R(ii)\; Each element of $\De$ can be written as
 $\sum_{s\in S}a_s\al_s$ with either $a_s\geq0$ for all $s\in S$
 or $a_s\leq0$ for all $s\in S$, but not in both ways.
 Accordingly, we write $\al>0$ or $\al<0$, and put
 $\De_+=\{\,\al\in\De\,|\,\al>0\,\}$ and
 $\De_-=\{\,\al\in\De\,|\,\al<0\,\}$.

 R(iii)\; For each $s\in S$, $s(\al_s)=-\al_s$ and
 $s(\De_+\setminus\{\al_s\})=\De_+\setminus\{\al_s\}$.

 R(iv)\; If $w\in W$ and $s,\,s'\in S$ satisfy
 $w(\al_{s'})=\al_s$, then $ws'w^{-1}=s$.

\begin{defn} 
 Define subsets $\re$, $\im$, $\re_\pm$, and $\im_\pm$ of $\De$
 by setting
\begin{align*}
 \re &:= \{\,w(\al_s)\;|\;w\in W,\,s\in S\,\},\qquad %
 \im := \De\setminus\re, \\
 \re_\pm &:= \re\cap\De_\pm,\qquad \im_\pm:=\im\cap\De_\pm.
\end{align*}
 Note that $W$ stabilizes $\re$ and $\pim$.
 We also set
\[ \vPhi(y):=\{\,\be\in\De_+\;|\;y^{-1}(\be)<0\,\} \]
 for each $y\in W$. Note that $\vPhi(y)\subset\pre$.
\end{defn}

\begin{thm}[\cite{vD2}] 
 The pair $(W,S)$ is a Coxeter system, i.e.,
 it satisfies the exchange condition.
 Moreover, if $y=s_1s_2\cdots s_n$ is a reduced expression of
 an element $y\in W\setminus\{1\}$, then
\[ \vPhi(y)=\{\,\al_{s_1},\,{s_1}(\al_{s_2}), %
 \,\dots,\,s_1\cdots s_{n-1}(\al_{s_n})\,\} \]
 and the elements of $\vPhi(y)$ displayed above are all distinct.
 In particular, $\sharp{\vPhi(y)}=\ell(y)$ for each $y\in W$,
 where $\ell\colon W\to{\bd Z}_{\geq0}$ is the length function.
\end{thm}

\noindent
{\it Remarks\/.}
 (1) For each Coxeter system $(W,S)$, a triplet $(V,\De,\Pi)$
 is called a {\it root system\/} of $(W,S)$ if
 it satisfies the conditions R(i)--R(iv).

 (2) Let $\sigma\colon W\to\mathrm{GL}(V)$ be
 the geometric representation of a Coxeter system $(W,S)$
 (cf. \cite{nB}), where $V$ is a real vector space having
 a basis $\Pi=\{\,\al_s\,|\,s\in S\,\}$.
 Then $(V,\De,\Pi)$ is a root system of $(W,S)$
 (cf. \cite{vD2}), where
 $\De=\{\,\sigma(w)(\al_s)\,|\,w\in W,\,s\in S\,\}$.
 We call it the root system associated with
 the geometric representation.

 (3) Let ${\fr g}$ be a Kac-Moody Lie algebra over ${\bd R}$
 with a Cartan subalgebra ${\fr h}$,
 $\De\subset{\fr h}^*\setminus\{0\}$ the root system of ${\fr g}$,
 $\Pi=\{\,\al_i\,|\,i\in{\bd I}\,\}$ a root basis of $\De$, and
 $W=\la{\,s_i\,|\,i\in{\bd I}\,}\ra\subset{\mathrm GL}({\fr h}^*)$
 the Weyl group of ${\fr g}$, where ${\fr h}^*$ is the dual space of
 ${\fr h}$ and $s_i$ is the simple reflection associated with $\al_i$
 (cf. \cite{vK}).
 Then $({\fr h}^*,\De,\Pi)$ is a root system of
 a Coxeter system $(W,S)$, where $S=\{\,s_i\,|\,i\in{\bd I}\,\}$.

\begin{defn} 
 For subsets $A,\,B\subset\De_+$ satisfying $B\subset A$,
 we say that $B$ is {\it convex in $A$\/}
 if it satisfies the following condition:

 $\mbox{C(i)}_{\ti{A}}\quad %
 \be,\,\ga\in B,\;\; \be+\ga\in A \;\Lra\; \be+\ga\in B$.

\noindent
 We also say that $B$ is {\it coconvex in $A$\/}
 if it satisfies the following condition:

 $\mbox{C(ii)}_{\ti{A}}\quad %
 \be,\,\ga\in A\setminus B,\;\; %
 \be+\ga\in A \;\Lra\; \be+\ga\in A\setminus B$.

\noindent
 Further we say that $B$ is {\it biconvex in $A$\/}
 if $B$ is both convex and coconvex in $A$.
 Note that $B$ is coconvex in $A$ if and only if
 $A\setminus B$ is convex in $A$ and that
 the condition $\mbox{C(ii)}_{\ti{A}}$
 is equivalent to the following condition:
\[ \be,\,\ga\in A,\;\; \be+\ga\in B %
 \;\Lra\; \be\in B \;\mbox{ or }\; \ga\in B. \]

 We simply say that
 $B$ is {\it convex\/}, {\it coconvex\/}, or {\it biconvex\/}
 if $B$ is convex in $\De_+$, coconvex in $\De_+$,
 or biconvex in $\De_+$, respectively.
 We also write C(i) and C(ii) instead of
 $\mbox{C(i)}_{\ti{\De_+}}$ and $\mbox{C(ii)}_{\ti{\De_+}}$,
 respectively.
 Further we call a subset $B\subset\pre$ a {\it real convex set\/},
 a {\it real coconvex set\/}, or a {\it real biconvex set\/}
 if $B$ is convex, coconvex, or biconvex, respectively.
 We denote by ${\fr B}$ the set of all finite biconvex sets, and
 by ${\fr B}^\infty$ the set of all infinite real biconvex sets.
\end{defn}

\begin{thm}[\cite{pP1}] 
 The assignment $y\,\mapsto\,\vPhi(y)$ defines
 an injective mapping between $W$ and ${\fr B}$.
 Moreover, if the root system $(V,\De,\Pi)$ satisfies
 the following two conditions then $\vPhi$ is surjective{\em:}

 {\em R(v)}\; each $\al\in\De_+\setminus\Pi$ can be written as
 $\be+\ga$ with $\be,\,\ga\in\De_+;$

 {\em R(vi)}\; there exists a mapping
 ${\mathrm ht}\colon\De_+\to{\bd R}_{>0}$ such that
 ${\mathrm ht}(\be+\ga)={\mathrm ht}(\be)+{\mathrm ht}(\ga)$
 if $\be+\ga\in\De_+$ for some $\be,\,\ga\in\De_+$.
\end{thm}

\noindent
{\it Remarks\/.}
 (1) Suppose that the root system $(V,\De,\Pi)$ satisfies
 the following two conditions instead of R(v) and R(vi):

 $\mbox{R(v)}'$\; each $\al\in\De_+\setminus\Pi$ can be written as
 $b\be+c\ga$ with $b,\,c\in{\bd R}_{\geq1}$ and $\be,\,\ga\in\De_+$;

 $\mbox{R(vi)}'$\; there exists a mapping
 ${\mathrm ht}\colon\De_+\to{\bd R}_{>0}$ such that
 ${\mathrm ht}(b\be+c\ga)= %
 b\,{\mathrm ht}(\be)+c\,{\mathrm ht}(\ga)$
 if $b\be+c\ga\in\De_+$ for some $b,\,c\in{\bd R}_{>0}$
 and $\be,\,\ga\in\De_+$.

\noindent
 Then $\vPhi$ is still surjective if ${\fr B}$ is changed into
 the set of all finite subsets $B\subset\De_+$ which satisfy
 the following two conditions:

 $\mbox{C(i)}'$\quad %
 $\be,\,\ga\in B,\;\; b,\,c\in{\bd R}_{>0},\;\; %
 b\be+c\ga\in\De_+ \;\Lra\; b\be+c\ga\in B$;

 $\mbox{C(ii)}'$\quad %
 $\be,\,\ga\in\De_+\setminus B,\;\; b,\,c\in{\bd R}_{>0},\;\; %
 b\be+c\ga\in\De_+ \;\Lra\; b\be+c\ga\in\De_+\setminus B$.

 (2) Let $(V,\De,\Pi)$ be the root system of
 a Coxeter system $(W,S)$ associated with
 the geometric representation.
 Then $(V,\De,\Pi)$ satisfies $\mbox{R(v)}'$ and $\mbox{R(vi)}'$.
 The condition $\mbox{R(v)}'$ is easily checked by
 reforming the proof of Proposition 2.1 in \cite{vD1}.
 Since $\Pi$ is a basis of $V$ we can define a mapping
 ${\mathrm ht}\colon\De_+\to{\bd R}_{>0}$ by setting
 ${\mathrm ht}(\al)=\sum_{s\in S}a_s$ for each $\al\in\De_+$,
 where $a_s$'s are unique non-negative real numbers such that
 $\al=\sum_{s\in S}a_s\al_s$.
 Then ${\mathrm ht}$ satisfies
 the required property in $\mbox{R(vi)}'$.

\begin{lem} 
 Let $B$ and $C$ be subsets of $\De_+$ satisfying
 $C=\vPhi(y)\amalg yB$ with $y\in W$.
 Then $B$ is biconvex if and only if $C$ is biconvex.
\end{lem}

\begin{pf}
 It is easy to see that
\begin{align}
 B = y^{-1}C\cap\De_+, \\ 
 \De_+\setminus B = \De_+\setminus y^{-1}C, \\ 
 y^{-1}(\De_+\setminus C) \subset \De_+\setminus B. 
\end{align}

 Suppose that $C$ is biconvex.
 To check the convexity of $B$, suppose that
 $\be\in B$ and $\ga\in B$ satisfy $\be+\ga\in\De_+$.
 Then $y(\be),\,y(\ga)\in C$.
 Since $y(\be+\ga)>0$ we have $y(\be+\ga)\in C\cap y\De_+$
 by the convexity of $C$, and hence $\be+\ga\in B$.
 To check the convexty of $\De_+\setminus B$, suppose that
 $\be\in\De_+$ and $\ga\in\De_+$ satisfy $\be+\ga\in B$.
 Then $y(\be+\ga)\in C$.
 In case $y(\be)>0$ and $y(\ga)>0$, we have either
 $y(\be)\in C\cap y\De_+$ or $y(\ga)\in C\cap y\De_+$
 by the convexity of $\De_+\setminus C$,
 and hence $\be\in B$ or $\ga\in B$.
 In case $y(\be)<0$ and $y(\ga)>0$,
 we have $-y(\be)\in\vPhi(y)\subset C$.
 Thus we get
 $y(\ga)=-y(\be)+\{y(\be+\ga)\}\in C\cap y\De_+$
 by the convexity of $C$, and hence $\ga\in B$.

 Suppose that $B$ is biconvex.
 To check the convexity of $C$, suppose that
 $\be\in C$ and $\ga\in C$ satisfy $\be+\ga\in\De_+$.
 By the convexity of $\vPhi(y)$ and $B$, we may assume that
 $\be\in\vPhi(y)$, $\ga\in yB$, and $\be+\ga\notin\vPhi(y)$.
 Since $-y^{-1}(\be)\in\De_+\setminus B$, the equality
 $y^{-1}(\ga)=y^{-1}(\be+\ga)-y^{-1}(\be)$ implies that
 $y^{-1}(\be+\ga)\in B$ by the convexity of $\De_+\setminus B$,
 and hence $\be+\ga\in yB\subset C$.
 To check the convexity of $\De_+\setminus C$, suppose that
 $\be\in\De_+\setminus C$ and $\ga\in\De_+\setminus C$
 satisfy $\be+\ga\in\De_+$.
 Then $y^{-1}(\be),\,y^{-1}(\ga)\in\De_+\setminus B$ by (2.3).
 Since $y^{-1}(\be+\ga)>0$ we have
 $y^{-1}(\be+\ga)\in\De_+\setminus B$
 by the convexity of $\De_+\setminus B$,
 and hence $\be+\ga\in\De_+\setminus C$ by (2.2).
\end{pf}

\begin{defn} 
 Let $\Lam$ be either the set ${\bd N}$ or the finite subset
 $\{\,1,\,\cdots,\,n\}\subset{\bd N}$ with $n\in{\bd N}$.
 For each sequence $\bs=(\bs(p))_{p\in\Lam}\in S^\Lam$,
 we define two mappings
 $z_{\bs}\colon\Lam\to W$ and $\phi_{\bs}\colon\Lam\to\re$
 by setting for each $p\in\Lam$
\[ z_{\bs}(p):=\bs(1)\cdots\bs(p),\qquad %
 \phi_{\bs}(p):=z_{\bs}(p-1)(\al_{\bs(p)}), \]
 where $z_{\bs}(0):=1$.
 We say that an element $\bs\in S^\Lam$ is {\it reduced\/}
 if $\ell(z_{\bs}(p))=p$ for all $p\in\Lam$,
 and denote by ${\cl W}^\Lam$ the subset of $S^\Lam$
 consisting of all reduced elements.

 In case $\Lam={\bd N}$, we write ${\cl W}^\infty$
 instead of ${\cl W}^\Lam$ and call an element of
 ${\cl W}^\infty$ a {\it infinite reduced word\/}.
 For a pair $(\bs,\,\bs')$ of infinite reduced words,
 we write $\bs\sim\bs'$ if there exists
 $(p_0\,,\,q_0)\in{\bd Z}_{\geq p}\times{\bd Z}_{\geq q}$
 for each $(p\,,\,q)\in{\bd N}^2$ such that
\[ \ell(z_{\bs}(p)^{-1}z_{\bs'}(p_0))=p_0-p,\quad %
 \ell(z_{\bs'}(q)^{-1}z_{\bs}(q_0))=q_0-q. \]
 Then $\sim$ is an equivalence relation on
 ${\cl W}^\infty$ (see \cite{kI}).
 Denote by $W^\infty$ the quotient set of ${\cl W}^\infty$
 relative to $\sim$, and by $[\bs]$ the coset including
 $\bs\in{\cl W}^\infty$.
 Define a mapping $\wtil{\vPhi}^\infty$ from $S^{\bd N}$ to
 the power set of $\pre$ by setting for each $\bs\in S^{\bd N}$
\[ \wtil{\vPhi}^\infty(\bs):= %
 \bigcup_{p\in{\bd N}}\vPhi(z_{\bs}(p)). \]

 In case $\Lam=\{\,1,\,\cdots,\,n\}$,
 we write ${\cl W}^n$ instead of ${\cl W}^\Lam$.
 We put ${\cl W}:=\amalg_{n\in{\bd N}}{\cl W}^n$ and call
 an element of ${\cl W}$ a {\it finite reduced word\/}.
 For each finite reduced word $\bs$, we put
 $[\bs]:=z_{\bs}(n)$ if $\bs\in{\cl W}^n$.
\end{defn}

\begin{prop}[\cite{kI}] 
 Let $\bs$ and $\bs'$ be elements of $S^{\bd N}$.

 {\em(1)} We have $\bs\in{\cl W}^\infty$ if and only if
 $\phi_{\bs}(p)>0$ for all $p\in{\bd N}$.

 {\em(2)} If $\bs\in{\cl W}^\infty$, then
 $\wtil{\vPhi}^\infty(\bs)=\{\,\phi_{\bs}(p)\;|\;p\in{\bd N}\,\}$
 and $\phi_{\bs}(p)$'s are all distinct.

 {\em(3)} If $\bs\in{\cl W}^\infty$, then
 $\wtil{\vPhi}^\infty(\bs)\in{\fr B}^\infty$.

 {\em(4)} If $\bs,\,\bs'\in{\cl W}^\infty$,
 then we have $\bs\sim\bs'$ if and only if
 $\wtil{\vPhi}^\infty(\bs)=\wtil{\vPhi}^\infty(\bs')$.
\end{prop}

\begin{defn} 
 Thanks to Proposition 2.7, we have an injective mapping
\[ \vPhi^\infty\colon W^\infty\,\lra\,{\fr B}^\infty,\quad %
 [\bs]\,\longmapsto\, %
 \vPhi^\infty([\bs]):=\wtil{\vPhi}^\infty(\bs). \]
\end{defn}

 We next define a left action of $W$ on $W^\infty$.

\begin{defn} 
 For each $x\in W$ and $\bs\in S^{\bd N}$, we set
\[ \wtil{\vPhi}^\infty(x,\bs):= %
 \{\, \be\in\pre \;|\; \exists p_0\in{\bd N};\, %
 \forall p\geq p_0,\,(xz_{\bs}(p))^{-1}(\be)<0 \,\}. \]
\end{defn}

\begin{lem}[\cite{kI}] 
 {\em(1)} If $\bs\in{\cl W}^\infty$, then
 $\wtil{\vPhi}^\infty(1,\bs)=\wtil{\vPhi}^\infty(\bs)$.

 {\em(2)} If $x\in W$ and $\bs\in{\cl W}^\infty$,
 then there exists an element $\bs'\in{\cl W}^\infty$ such that
 $\wtil{\vPhi}^\infty(\bs')=\wtil{\vPhi}^\infty(x,\bs)$.
 More precisely, a required $\bs'\in S^{\bd N}$ is constructed
 by applying the following procedure {\em Step 1--3}.

 {\em Step 1.} Choose a non-negative integer $p_0$ satisfying
\[ \vPhi(x^{-1})\cap{\wtil{\vPhi}^\infty(\bs)}\subset %
 {\vPhi(z_{\bs}(p_0))}. \]

 {\em Step 2.}
 In case $xz_{\bs}(p_0)=1$, put $\bs'(p):=\bs(p_0+p)$
 for each $p\in{\bd N}$.
 In case $xz_{\bs}(p_0)\neq1$, choose a reduced expression
 $xz_{\bs}(p_0)=\bs'(1)\cdots\bs'(l_0)$, and put
 $\bs'(l_0+p):=\bs(p_0+p)$ for each $p\in{\bd N}$.

 {\em Step 3.} Set $\bs'\!:=\!(\bs'(p))_{p\in{\bd N}}$.

 {\em(3)} Let $\bs$ and $\bs'$ be elements of ${\cl W}^\infty$.
 Then the following two properties hold{\em:}
\begin{itemize}
\item[(i)] $\wtil{\vPhi}^\infty(\bs)=\wtil{\vPhi}^\infty(\bs') %
 \,\Lra\,\wtil{\vPhi}^\infty(x,\bs)=\wtil{\vPhi}^\infty(x,\bs');$
\item[(ii)] $\wtil{\vPhi}^\infty(y,\bs)=\wtil{\vPhi}^\infty(\bs') %
 \,\Lra\,\wtil{\vPhi}^\infty(xy,\bs)=\wtil{\vPhi}^\infty(x,\bs')$
\end{itemize}
 for each pair of elements $x,\,y$ of $W$.
\end{lem}

\begin{defn} 
 By Proposition 2.7(4) and Lemma 2.10, we have
 a left action of $W$ on $W^\infty$ satisfying $x.[\bs]=[\bs']$
 if $x\in W$ and $\bs,\,\bs'\in{\cl W}^\infty$ satisfy
 $\wtil{\vPhi}^\infty(x,\bs)=\wtil{\vPhi}^\infty(\bs')$.
\end{defn}

\begin{prop}[\cite{kI}] 
 If $x\in W$ and $\bs\in{\cl W}^\infty$ then
\[ \vPhi^\infty(x.[\bs])= %
 \{\vPhi(x)\setminus(-\Omega)\}\amalg %
 \{x\vPhi^\infty([\bs])\setminus\Omega\}, \]
 where $\Omega:=x\vPhi^\infty([\bs])\cap\nre$.
 In particular, if $x\vPhi^\infty([\bs])\subset\pre$ then
\[ \vPhi^\infty(x.[\bs])= %
 \vPhi(x)\amalg{x\vPhi^\infty([\bs])}. \]
\end{prop}

%% file: sec3.tex
\section{Definition of convex orders}
 Before defining convex orders, we prepare notation for
 total orders on sets. Let $\preq$ be a total order on a set $R$.
 Denote by $\preq^{op}$ the total order on $R$ satisfying
 $r\preq^{op}r'$ if and only if $r'\preq r$ for $r,\,r'\in R$,
 and call it {\it the opposite of $\preq$\/}.
 For each subset $A\subset R$, we denote by $\preq_{\ti{A}}$
 the total order on $A$ satisfying $a\preq_{\ti{A}}a'$
 if and only if $a\preq a'$ for $a,\,a'\in A$, and call it
 {\it the restriction of $\preq$ to $A$\/}.
 Conversely, if the restriction $\preq_{\ti{A}}$ is
 equal to a total order $\preq'$ on $A$,
 we call $\preq$ {\it an extension of $\preq'$ to $R$\/}.
 For $r,\,r'\in R$, we write $r\pr r'$
 if $r\preq r'$ and $r\neq r'$.
 For subsets $A,\,B\subset R$, we write $A\pr B$
 if $a\pr b$ for each $(a,b)\in A\times B$.
 Recall that $\preq$ is called a {\it well-order\/}
 if each non-empty subset $T\subset R$ has
 the minimum element relative to $\preq_{\ti{T}}$.

\begin{defn} 
 Let $\preq$ be a total order on a set $R$.

 (1) If $R$ is finite and $\sharp R=n$, we denote by
 $\phi_\preq$ the unique order isomorphism from
 the usually ordered set $\{\,1,\,\cdots,\,n\,\}$ to $R$.

 (2) We say that $\preq$ is {\it of 1-row type\/} if $R$ is order
 isomorphic to the usually ordered set ${\bd N}$, and denote by
 $\phi_\preq$ the unique order isomorphism from ${\bd N}$ to $R$.

 (3) We say that $\preq$ is {\it of n-row type\/} if there exists
 a unique decomposition $R=\amalg_{i=1}^nR_i\ft{(\preq)}$ with
 a unique positive integer $n\in{\bd N}$ such that
\begin{itemize}
 \item[(i)] the restriction of $\preq$ to $R_i\ft{(\preq)}$
 is of 1-row type for each $i=1,\dots,n$,
 \item[(ii)] $R_i\ft{(\preq)}\pr R_{i'}\ft{(\preq)}$
 if and only if $i<i'$ for $i,\,i'=1,\dots,n$.
\end{itemize}
\end{defn}

\begin{lem} 
 {\em(1)} Let $\preq$ be a total order on an infinite set $R$.
 Then $\preq$ is of 1-row type if and only if
 $R_\xi:=\{\,r\in R\,|\,r\preq\xi\,\}$
 is finite for each $\xi\in R$.

 {\em(2)} Let $(R,\preq)$ be a well-ordered non-empty set
 without the maximum element.
 Then there exists a unique decomposition
 $R=\amalg_{\lam\in\Lam}R_\lam\ft{(\preq)}$ with a uniquely
 determined {\em(}up to order isomorphism\/{\em)} well-ordered
 non-empty set $(\Lam,\wt{\preq})$ such that
\begin{itemize}
 \item[(i)] the restriction of $\preq$ to $R_\lam\ft{(\preq)}$
 is of 1-row type for each $\lam\in\Lam$,
 \item[(ii)] $R_\lam\ft{(\preq)}\pr R_{\lam'}\ft{(\preq)}$
 if and only if $\lam\,\wt{\pr}\,\lam'$ for $\lam,\,\lam'\in\Lam$.
\end{itemize}
\end{lem}

\begin{pf}
 (1) Clear.

 (2) This follows from (1) and the standard results for
 well-ordered sets.
\end{pf}

\begin{defn} 
 Let $(V,\De,\Pi)$ be a root system of
 a Coxeter system $(W,S)$, and $\preq$ a total order on
 a subset $B\subset\De_+$.
 We say that $\preq$ is a {\it convex order\/}
 if it satisfies the following conditions:

 $\mbox{CO(i)}\quad %
 (\be,\,\ga)\in B^2\setminus(\pim)^2,\;\; %
 \be\,\pr\,\ga,\;\; \be+\ga\in B \;\Lra\; %
 \be\,\pr\,\be+\ga\,\pr\,\ga$;

 $\mbox{CO(ii)}\quad %
 \be\in B,\;\; \ga\in\De_+\setminus B,\;\; %
 \be+\ga\in B \;\Lra\; \be\,\pr\,\be+\ga$.

\noindent
 We also say that $\preq$ is an {\it opposite convex order\/}
 if the opposite $\preq^{op}$ of $\preq$ is a convex order.
 In addition, we denote by $\mbox{CO(ii)}_A$
 the following condition concerning $\preq$:
\[ \be\in B,\;\; \ga\in A\setminus B,\;\; \be+\ga\in B %
 \;\Lra\; \be\,\pr\,\be+\ga, \]
 where $A$ is a subset of $\De_+$. Note that
 $\mbox{CO(ii)}_{\De_+}$ is the same condition as CO(ii).
\end{defn}

\noindent
{\it Remark\/.}
 In case $\De_+$ is the positive root system relative to
 a root basis of a Kac-Moody Lie algebra, each root
 $\de\in\pim$ satisfies $n\de\in\pim$ for all $n\in{\bd N}$
 (cf. \cite{vK}). This fact denies the existence of
 a total order $\preq$ on $\pim$ satisfying $\be\pr\be+\ga\pr\ga$
 for each pair $(\be,\ga)\in(\pim)^2$ such that
 $\be\pr\ga$ and $\be+\ga\in\pim$.
 This is the reason why we subtract $(\pim)^2$ from $B^2$
 in the condition CO(i).

\noindent
{\it Example\/.}
 In case $\De$ is the root system of the affine Lie algebra of
 type $A_2^{(1)}$, the following infinite subset of $\pre$ is
 a maximal real biconvex set:
\[ \De\ft{(1,-)}= %
 \{\,m\de-\al_1,\,m\de-\al_1-\al_2,\,m\de-\al_2\;|\;m\in{\bd N}\,\}, \]
 where we use notation in Section 4.

 The following total order $\preq$ on $\De\ft{(1,-)}$ is
 a 1-row type convex order:
\begin{align*}
  \de-\al_1-\al_2 &\pr  \de-\al_1 \pr %
 2\de-\al_1-\al_2  \pr  \de-\al_2 \pr \\
 3\de-\al_1-\al_2 &\pr 2\de-\al_1 \pr %
 4\de-\al_1-\al_2  \pr 2\de-\al_2 \pr \\
 5\de-\al_1-\al_2 &\pr 3\de-\al_1 \pr %
 6\de-\al_1-\al_2  \pr 3\de-\al_2 \pr\cdots\cdot\cdot.
\end{align*}

 The following total order $\preq$ on $\De\ft{(1,-)}$ is
 a 2-row type convex order:
\begin{align*}
  \de-\al_1-\al_2 &\pr  \de-\al_1 \pr %
 2\de-\al_1-\al_2  \pr 2\de-\al_1 \pr %
 3\de-\al_1-\al_2  \pr 3\de-\al_1 \pr\cdots\cdot\cdot \\
 \pr \de-\al_2 &\pr 2\de-\al_2 \pr 3\de-\al_2 \pr\cdots\cdot\cdot
\end{align*}
 such that
\begin{align*}
 R_1\ft{(\preq)} &= %
 \{\,m\de-\al_1,\,m\de-\al_1-\al_2\;|\;m\in{\bd N}\,\}, \\
 R_2\ft{(\preq)} &= \{\,m\de-\al_2\;|\;m\in{\bd N}\,\}.
\end{align*}

\begin{lem} 
 {\em(1)} An opposite convex order on $\De_+$
 is a convex order on $\De_+$.

 {\em(2)} Let $\preq$ be a convex order on a set $C$, and
 $B$ is a subset of $C$ such that $B\pr C\setminus B$.
 Then the restriction $\preq_{\ti{B}}$ is a convex order on $B$.
 Moreover, if $C$ is biconvex, then $B$ is biconvex.
\end{lem}

\begin{pf}
 (1)\; Clear.

 (2)\; It is clear that
 $\preq_{\ti{B}}$ satisfies the condition CO(i).
 To check the property CO(ii) of $\preq_{\ti{B}}$, suppose that
 $\be\in B$ and $\ga\in\De_+\setminus B$ satisfy $\be+\ga\in B$.
 In case $\ga\in\De_+\setminus C$, we have
 $\be\pr_{\ti{B}}\be+\ga$ by the property CO(ii) of $\preq$.
 In case $\ga\in C\setminus B$, we have $\be\pr\ga$ since
 $B\preq C\setminus B$, and hence $\be\pr_{\ti{B}}\be+\ga\pr\ga$
 by the property CO(i) of $\preq$.

 In addition, we assume that $C$ is biconvex.
 To check the convexity of $B$, suppose that
 $\be\in B$ and $\ga\in B$ satisfy $\be+\ga\in\De_+$ and $\be\pr\ga$.
 We see that $\be+\ga\in C$ by the convexity of $C$ and
 that $\be\pr\be+\ga\pr\ga$ by the property CO(i) of $\preq$.
 Thus $\be+\ga\in B$ since $B\pr C\setminus B$.
 To check the convexity of $\De_+\setminus B$, suppose that
 $\be\in\De_+$ and $\ga\in\De_+$ satisfy $\be+\ga\in B$.
 By the convexity of $\De_+\setminus C$,
 we may assume that $\be\in C$.
 In case $\ga\in C$, we have
 either $\be\pr\be+\ga\pr\ga$ or $\ga\pr\be+\ga\pr\be$
 by the property CO(i) of $\preq$,
 and hence $\be\in B$ or $\ga\in B$
 since $B\pr C\setminus B$.
 In case $\ga\in\De_+\setminus C$, we have
 $\be\pr\be+\ga$ by the property CO(ii) of $\preq$,
 and hence $\be\in B$ since $B\pr C\setminus B$.
\end{pf}

\begin{prop} 
 Suppose that the root system $(V,\De,\Pi)$ satisfies
 the conditions {\em R(v)} and {\em R(vi) (see Theorem 2.4)}.

 {\em(1)(\cite{pP1})} Suppose that $B$ is a non-empty finite
 biconvex set such that $B=\vPhi([\bs])$ for some $\bs\in{\cl W}$.
 Define a total order $\preq$ on $B$ by setting
 $\phi_{\bs}(p)\preq\phi_{\bs}(q)$ for each $p\leq q$.
 Then $\preq$ is a convex order on $B$ such that
 $\phi_{\ti{\preq}}=\phi_{\bs}$.
 Moreover, each convex order $\preq$ on a non-empty finite
 biconvex set $C$, there exists a unique element
 $\bs\in{\cl W}$ such that $C=\vPhi([\bs])$ and
 $\phi_{\ti{\preq}}=\phi_{\bs}$.

 {\em(2)} Suppose that $B$ is an infinite real biconvex set
 such that $B=\vPhi^\infty([\bs])$ for some $\bs\in{\cl W}^\infty$.
 Define a total order $\preq$ on $B$ by setting
 $\phi_{\bs}(p)\preq\phi_{\bs}(q)$ for each $p\leq q$.
 Then $\preq$ is a 1-row type convex order on $B$ such that
 $\phi_{\ti{\preq}}=\phi_{\bs}$.
 Moreover, for each 1-row type convex order $\preq$ on an
 infinite real biconvex set $C$, there exists a unique element
 $\bs\in{\cl W}^\infty$ such that $C=\vPhi^\infty([\bs])$ and
 $\phi_{\ti{\preq}}=\phi_{\bs}$.
\end{prop}

\begin{pf}
 (1)\; See \cite{pP1}.

 (2)\; The proof of the first assertion is almost similar to
 that of the first assertion in (1). See \cite{pP1}.
 Let us prove the second assertion. For each $p\in{\bd N}$, put
 $C_p:=\{\,\phi_{\ti{\preq}}(q)\,|\,1\leq q\leq p\,\}$ and
 denote by $\preq_p$ the restriction of $\preq$ to $C_p$.
 Then we see that $C_p$ is a non-empty finite biconvex set and
 $\preq_p$ is a convex order on $C_p$ by Lemma 3.4(2).
 Hence, by (1), there exists a unique element
 $\bs_p\in{\cl W}$ such that $C_p=\vPhi([\bs_p])$ and
 $\phi_{\ti{\preq_p}}=\phi_{\bs_p}$.
 Put $\bs(p):=\bs_p(p)$ for each $p\in{\bd N}$ and
 $\bs:=(\bs(p))_{p\in{\bd N}}\in S^{\bd N}$.
 Then $\bs$ is the unique element of ${\cl W}^\infty$ such that
 $C=\vPhi^\infty([\bs])$ and $\phi_{\ti{\preq}}=\phi_{\bs}$.
\end{pf}

\noindent
{\it Remark\/.}
 In case the root system $(V,\De,\Pi)$ satisfies the conditions
 $\mbox{R(v)}'$ and $\mbox{R(vi)}'$ instead of R(v) and R(vi),
 then the parts (1) and (2) of the above proposition are still valid
 if the properties C(i) and C(ii) of $B$ are changed into
 the conditions $\mbox{C(i)}'$ and $\mbox{C(ii)}'$
 (see remarks below Theorem 2.4).

\begin{prop} 
 Let $B$ and $C$ be subsets of $\De_+$ satisfying
 $C=\vPhi(y)\amalg yB$ with $y\in W$.

 {\em(1)} If there exist a convex order $\preq_1$ on $\vPhi(y)$
 and a convex order $\preq'$ on $B$, then there exists a unique
 convex order $\preq$ on $C$ such that
\begin{itemize}
 \item[(i)] $\vPhi(y)\pr yB$,
 \item[(ii)] the restriction of $\preq$ to $\vPhi(y)$ is $\preq_1$,
 \item[(iii)] $y(\xi)\preq y(\eta)$ if and only if
 $\xi\preq'\eta$ for $\xi,\,\eta\in B$.
\end{itemize}

 {\em(2)} If there exists a convex order $\preq$ on $C$
 with the above property {\em(i)}, then
 there exists a unique convex order $\preq'$ on $B$
 with the above property {\em(iii)}.
\end{prop}

\begin{pf}
 (1)\; Since the properties (i)--(iii)
 determine a unique total order $\preq$ on $C$,
 it suffices to check that $\preq$ is a convex order on $C$.
 To check the property CO(i) of $\preq$, suppose that
 $\be\in C$ and $\ga\in C$ satisfy $\be\pr\ga$ and $\be+\ga\in C$.
 By the properties (ii) and (iii),
 we may assume that $\be\in\vPhi(y)$ and $\ga\in yB$.
 In case $\be+\ga\in\vPhi(y)$, we have
 $\be+\ga\pr\ga$ by the property (i).
 Moreover we have $\be\preq_1\be+\ga$
 by the property CO(ii) of $\preq_1$,
 and hence $\be\pr\be+\ga\pr\ga$.
 In case $\be+\ga\in yB$, we have
 $\be\pr\be+\ga$ by the property (i).
 Since $-y^{-1}(\be)\in\De_+\setminus B$, the equality
 $y^{-1}(\ga)=y^{-1}(\be+\ga)-y^{-1}(\be)$ implies that
 $y^{-1}(\be+\ga)\pr'y^{-1}(\ga)$
 by the property CO(ii) of $\preq'$,
 and hence $\be+\ga\pr\ga$ by the property (iii).
 Thus we get $\be\pr\be+\ga\pr\ga$.
 We next check the property CO(ii) of $\preq$.
 Suppose that $\be\in C$ and $\ga\in\De_+\setminus C$
 satisfy $\be+\ga\in C$.
 By the property CO(ii) of $\preq_1$,
 we may assume that $\be\in yB$.
 By (2.3), we have $y^{-1}(\ga)\in\De_+\setminus B$,
 and hence $\be+\ga\in yB$ since $y^{-1}(\be+\ga)>0$.
 Thus we get $y^{-1}(\be)\pr'y^{-1}(\be+\ga)$
 by the property CO(ii) of $\preq'$, and hence
 $\be\pr\be+\ga$ by the property (iii).

 (2)\; It is clear that there exists a unique total order
 $\preq'$ on $B$ with the property (iii).
 Let us check that $\preq'$ is a convex order on $B$.
 The property CO(i) is clear.
 To check the property CO(ii), suppose that
 $\xi\in B$ and $\eta\in\De_+\setminus B$ satisfy $\xi+\eta\in B$.
 In case $y(\eta)>0$, we have
 $y(\eta)\in\De_+\setminus C$ by (2.2).
 Since $y(\xi),\,y(\xi+\eta)\in C$ we have
 $y(\xi)\pr y(\xi+\eta)$ by the property CO(ii) of $\preq$,
 and hence $\xi\pr'\xi+\eta$.
 In case $y(\eta)<0$, we have $-y(\eta)\in\vPhi(y)$,
 and hence $-y(\eta)\pr y(\xi+\eta)$ by the property (i).
 Thus we get $-y(\eta)\pr y(\xi)\pr y(\xi+\eta)$
 by the property CO(i) of $\preq$,
 and hence $\xi\pr'\xi+\eta$.
\end{pf}

%% file: sec4.tex
\section{Notation for untwisted affine cases}
 In this section, we prepare notation for untwisted affine cases.
 Let $l$ be a positive integer, and put
 ${\bd I}=\{0,1,\dots,l\}$ and ${\SI}=\{1,\dots,l\}$.
 Let ${\mathrm A}=(a_{ij})_{i,j\in{\bd I}}$ be a generalized
 Cartan matrix of affine type $X_l^{(1)}$ such that
 $(a_{ij})_{i,j\in{\SI}}$ is the Cartan matrix of type $X_l$,
 where $X=A,B,\dots,\mbox{ or }G$.
 Let $({\fr h},\Pi,\Pi^\vee)$ be a minimal realization of
 ${\mathrm A}$ over ${\bd R}$, that is, a triplet consisting of
 a $(l+2)$-dimensional real vector space ${\fr h}$ and
 linearly independent subsets
 $\Pi=\{\al_i\,|\,i\in{\bd I}\}\subset{\fr h}^*$,
 $\Pi^\vee=\{\al_i^\vee\,|\,i\in{\bd I}\}\subset{\fr h}$
 satisfying $\la{\,\al_i\,,\,\al_j^\vee\,}\ra=a_{ij}$
 for each $i,\,j\in{\bd I}$.
 Let ${\fr g}$ be the Kac-Moody Lie algebra associated with
 $({\fr h},\Pi,\Pi^\vee)$ or the affine Lie algebra
 of type $X_l^{(1)}$, $\De\subset{\fr h}^*\setminus\{0\}$
 the set of all roots of ${\fr g}$,
 $\re$ (resp. $\im$) the set of all real (resp. imaginary) roots,
 and $W=\la{\,s_i\,|\,i\in{\bd I}\,}\ra$ the Weyl group
 of ${\fr g}$, where $s_i\in\mathrm{GL}({\fr h}^*)$
 is the reflection with respect to $\al_i$.
 Further, let $\De_+$ (resp. $\De_-$) be the set of
 all positive (resp. negative) roots relative to $\Pi$,
 and ${\mathrm ht}\colon\De_+\to{\bd N}$
 the height function on $\De_+$.

 Set ${\sPi}:=\{\al_i\,|\,i\in{\SI}\}$,
 ${\st{\circ}{\fr h}}{}^*:=\oplus_{i\in{\SI}}{\bd R}\al_i$,
 ${\sW}:=\la{\,s_i\,|\,i\in{\SI}\,}\ra$,
 ${\sDe}:={\sW}(\st{\circ}{\Pi})$,
 and ${\sDe}_\pm:={\sDe}\cap\De_\pm$.
 Note that ${\sDe}$ is a root system of type $X_l$ with
 a root basis ${\sPi}$ and the Weyl group ${\sW}$.

 Let $a_i$ ($i\in{\bd I}$) be the labels of ${\mathrm A}$,
 and put $\de=\sum_{i\in{\bd I}}a_i\al_i$.
 Then $\la{\de\,,\,\al_i^\vee}\ra=0$ for all $i\in{\bd I}$,
 and hence $w(\de)=\de$ for all $w\in W$.
 By the assumption for ${\mathrm A}$, we have $a_0=1$ and
 $\de=\al_0+\theta$, where $\theta=\sum_{i\in{\SI}}a_i\al_i$
 is the highest root of ${\sDe}$. Moreover, we have
\[ \re=\{\,m\de+\vep\;|\;m\in{\bd Z},\,\vep\in{\sDe}\,\},\quad %
 \im=\{\,m\de\;|\;m\in{\bd Z}\setminus\{0\}\,\}. \]

 Let $(\cdot\,|\,\cdot)$ be the standard symmetric bilinear form
 on ${\fr h}^*$, scaled so that $(\theta\,|\,\theta)=2$.
 Note that $(\de\,|\,\al_i)=0$ for all $i\in{\bd I}$ and that
 the restriction of the form $(\cdot\,|\,\cdot)$
 to ${\st{\circ}{\fr h}}{}^*$ is positive-definite.
 For each $\lam\in{\fr h}^*$, we denote by $\ol{\lam}$
 the image of $\lam$ by the orthogonal projection onto
 ${\st{\circ}{\fr h}}{}^*$.
 Note that each $\be\in\De$ can be uniquely written as
 $m\de+\ol{\be}$ with $m\in{\bd Z}$ and
 $\ol{\be}\in{\sDe}\amalg\{0\}$.

 For each $\al\in\re$, we denote by $s_\al$ the reflection
 with respect to $\al$. For each $\lam\in{\fr h}^*$,
 we define an element $t_\lam\in\mathrm{GL}({\fr h}^*)$ by setting
\[ t_\lam(\mu)=\mu+(\mu\,|\,\de)\lam-\{(\mu\,|\,\lam)+ %
 \frac{1}{2}(\lam\,|\,\lam)(\mu\,|\,\de)\}\de \]
 for each $\mu\in{\fr h}^*$.
 We have $t_\lam(\mu)=\mu-(\mu\,|\,\lam)\de$
 for each $\mu\in{\fr h}_0^*$, where
 ${\fr h}_0^*:=\oplus_{i\in{\bd I}}{\bd R}\al_i$.

 For an arbitrary subset ${\bd J}\subset{\SI}$, we set
\begin{align*}
 {\sPi}_{\tiJ} &:= \{\,\al_j\;|\;j\in{\bd J}\,\}\subset{\sPi}, %
 \qquad{\sW}_{\tiJ} := \la{\,r_j\;|\;j\in{\bd J}\,}\ra\subset{\sW}, \\
 {\sDe}_{\tiJ} &:= {\sW}_{\tiJ}({\sPi}_{\tiJ})\subset{\sDe},\qquad %
 {\sDe}_{\tiJ\pm} := {\sDe}_{\tiJ}\cap{\sDe}_\pm.
\end{align*}
 Note that ${\sDe}_{\tiJ}$ is a semi-simple root system with
 a root basis ${\sPi}_{\tiJ}$ and the Weyl group ${\sW}_{\tiJ}$
 if ${\bd J}\neq\emptyset$.

\begin{defn} 
 For each subset ${\bd J}\subset{\SI}$, we set
\begin{align*}
 \re_{\tiJ} &:= \{\,\be\in\re\;|\;\ol{\be}\in{\sDe}_{\tiJ}\,\}, %
 \qquad\De_{\tiJ} := \re_{\tiJ}\amalg\im, \\
 \re_{\tiJ\pm} &:= \re_{\tiJ}\cap\De_\pm,\qquad %
 {\De}_{\tiJ\pm} := {\De}_{\tiJ}\cap\De_\pm.
\end{align*}
 Note that $\re_{\tiJ}$ and $\De_{\tiJ}$ are symmetric.
 For each non-empty subset ${\bd J}\subset{\SI}$, let
\[ {\sDe}_{\tiJ} = %
 \coprod_{c=1}^{\ti{{\mathrm C}({\bd J})}}{\sDe}_{{\tiJ}_c} \]
 be the irreducible decomposition of ${\sDe}_{\tiJ}$ with
 a unique positive integer $\ft{{\mathrm C}({\bd J})}$, and
 $\theta_{{\tiJ}_c}$ the highest root of ${\sDe}_{{\tiJ}_c}$
 relative to the basis ${\sPi}_{{\tiJ}_c}$ for each
 $c=1,\dots,\ft{{\mathrm C}({\bd J})}$. We set
\begin{align*}
 \Pi_{\tiJ_c} &:= {\sPi}_{\tiJ_c}\amalg\{\de-\theta_{\tiJ_c}\} %
 \,\mbox{ for }\, c=1,\dots,\ft{{\mathrm C}({\bd J})}, \\
 \Pi_{\tiJ} &:= \coprod_{c=1}^{\ti{{\mathrm C}({\bd J})}} %
 {\Pi}_{\tiJ_c},\qquad S_{\tiJ} := \{\,s_\al\;|\;\al\in\Pi_{\tiJ}\,\}.
\end{align*}
 For each $s\in S_{\tiJ}$, we denote by $\al_s$
 the unique element of $\Pi_{\tiJ}$ such that $s=s_{\al_s}$.
 Note that
 $\Pi_{\tiJ}=\{\,\al_s\,|\,s\in S_{\tiJ}\,\}\subset\re_{\tiJ+}$
 and $S_{\tiJ}\subset W$.
 We denote by $W_{\tiJ}$ the subgroup of $W$ generated by
 $S_{\tiJ}$, and by $V_{\tiJ}$ the subspace of ${\fr h}^*$
 spanned by $\Pi_{\tiJ}$. Note that
\begin{equation*}
 \re_{\tiJ}=\coprod_{c=1}^{\ti{{\mathrm C}({\bd J})}}\re_{\tiJ_c},
 \qquad W_{\tiJ}=\prod_{c=1}^{\ti{{\mathrm C}({\bd J})}}W_{\tiJ_c} %
 \;\mbox{(direct product)}.
\end{equation*}
 We also put $W_\emptyset:=\{1\}\subset W$.
\end{defn}

\begin{prop}[\cite{kI}] 
 For each non-empty subset ${\bd J}\subset{\SI}$,
 the pair $(W_{\tiJ},S_{\tiJ})$ is a Coxeter system and
 the triplet $(V_{\tiJ},\De_{\tiJ},\Pi_{\tiJ})$ is
 a root system of $(W_{\tiJ},S_{\tiJ})$ with
 the properties {\em R(v)} and {\em R(vi)}.
\end{prop}

\begin{defn} 
 For each subset ${\bd J}\subset{\SI}$, we define
 two subsets of the power set of $\De_{\tiJ+}$ by setting
\begin{align*}
 {\fr B}_{\tiJ} &:= \{\,B\subset\De_{\tiJ+} %
 \;|\; \mbox{$B$ is finite and biconvex in $\De_{\tiJ+}$}\,\}, \\
 {\fr B}_{\tiJ}^\infty &:= \{\,B\subset\re_{\tiJ+}\;|\; %
 \mbox{$B$ is infinite and biconvex in $\De_{\tiJ+}$}\,\}.
\end{align*}
 Note that ${\fr B}_\emptyset=\{\,\emptyset\,\}$ and
 ${\fr B}_\emptyset^\infty=\emptyset$.
 For each non-empty subset ${\bd J}\subset{\SI}$,
 we denote by ${\cl W}_{\tiJ}^\infty$ the set of
 all infinite reduced words of the Coxeter system
 $(W_{\tiJ},S_{\tiJ})$.
 Further, let $W_{\tiJ}^\infty$ be the quotient set of
 ${\cl W}_{\tiJ}^\infty$ obtained by applying Definition 2.6 for
 the Coxeter system $(W_{\tiJ},S_{\tiJ})$, and
 $\vPhi_{\tiJ}^\infty\colon %
 W_{\tiJ}^\infty\to{\fr B}_{\tiJ}^\infty$ the injective mapping
 obtained by applying Definition 2.8 for the root system
 $(V_{\tiJ},\De_{\tiJ},\Pi_{\tiJ})$ of $(W_{\tiJ},S_{\tiJ})$.
\end{defn}

\begin{cor}[\cite{CP}] 
 For each subset ${\bd J}\subset{\SI}$, the assignment
 $y\,\mapsto\,\vPhi_{\tiJ}(y):=\vPhi(y)\cap\De_{\tiJ+}$ defines
 a bijective mapping between $W_{\tiJ}$ and ${\fr B}_{\tiJ}$.
\end{cor}

\begin{pf}
 Since
 $\vPhi_{\tiJ}(y)=\{\,\be\in\De_{\tiJ+}\;|\;y^{-1}(\be)<0\,\}$,
 this follows immediately from Theorem 2.4 and Proposition 4.2.
\end{pf}

\noindent
{\it Remark\/.}
 The above corollary was stated by P.~Cellini and P.~Papi
 in the proof of Theorem 3.12 in \cite{CP}
 with an outline of the proof.

\begin{lem}[\cite{vK}] 
 Let ${\bd J}$ be a non-empty subset of ${\SI}$.
 Set $\check{\al}_j=\frac{2\al_j}{(\al_j\,|\,\al_j)} %
 \in{\st{\circ}{\fr h}}{}^\ast$ for each $j\in{\bd J}$,
 and put $M_{\tiJ}=\oplus_{j\in{\bd J}}{\bd Z}\check{\al}_j$
 and $T_{\tiJ}=\{\,t_{\lam}\,|\,\lam\in M_{\tiJ}\,\}$.
 Then $T_{\tiJ}$ is a normal subgroup of $W_{\tiJ}$
 such that $T_{\tiJ}\simeq M_{\tiJ}$ and
 $W_{\tiJ}=T_{\tiJ}\,{\sW}_{\tiJ}$.
\end{lem}

 For each $x\in W_{\tiJ}$, we denote by $\ol{x}$ the unique element
 of ${\sW}_{\tiJ}$ satisfying $x\in T_{\tiJ}\ol{x}$.
 Note that the mapping
 $\ol{\,\cdot\,}\colon W_{\tiJ}\to{\sW}_{\tiJ},\,x\mapsto\ol{x}$
 is a group homomorphism satisfying $\ol{x(\lam)}=\ol{x}(\ol{\lam})$
 and $\ol{s_\al}=s_{\ol{\al}}$ for all $x\in W_{\tiJ}$,
 $\lam\in V_{\tiJ}$, and $\al\in\re_{\tiJ}$.

%% file: sec5.tex
\section{Infinite real biconvex sets and infinite reduced words}
 In this section, we introduce a parameterization of
 the set ${\fr B}_{\tiJ}^\infty$ and a relationship between
 ${\fr B}_{\tiJ}^\infty$ and infinite reduced words of
 the Coxeter system $(W_{\tiJ},S_{\tiJ})$, where ${\bd J}$
 is an arbitrary non-empty subset of ${\SI}$.
 All statements in this section were already proved
 in case ${\bd J}={\SI}$ in \cite{kI}.
 We omit the proofs of all statements, since we can prove them
 for arbitrary ${\bd J}$ as in case ${\bd J}={\SI}$.

\begin{defn} 
 For each $\vep\in{\sDe}$ and $P\subset{\sDe}$,
 we define $\la{\vep}\ra,\,\la{P}\ra\subset\pre$ by setting
\[ \la{\vep}\ra:=\{\,m\de+\vep\;|\;m\in{\bd Z}_{\geq0}\,\}\cap\pre, %
 \qquad \la{P}\ra:=\coprod_{\vep\in P}\la{\vep}\ra. \]
\end{defn}

\begin{lem} 
 {\em(1)} For each $P\subset{\sDe}$ and $x\in W$,
 the following properties hold{\em:}
\[ {\mathrm(i)}\; %
 \la{P}\ra=\{\,\be\in\De_+\;|\;\ol{\be}\in{P}\,\},\quad %
 {\mathrm(ii)}\; \ol{x\la{P}\ra}\subset\ol{x}P,\quad %
 {\mathrm(iii)}\; x\la{P}\ra\doteq\la{\ol{x}P}\ra. \]

 {\em(2)} If a subset $B\subset\re_{\tiJ+}$ is biconvex in
 $\De_{\tiJ+}$ and a subset $P\subset{\sDe}_{\tiJ}$ satisfies
 $\la{P}\ra\,\dot{\subset}\,B$, then
 $\la{P}\ra\subset B$ and $\la{-P}\ra\cap B=\emptyset$.
\end{lem}

 For each ${\bd K}\subset{\bd J}$, we denote by
 ${\sW}{}^{\tiK}_{\tiJ}$ the minimal coset representatives for
 the set ${\sW}_{\tiJ}/{\sW}_{\tiK}$ of all right cosets.
 In case ${\bd J}={\SI}$ we denote it simply by ${\sW}{}^{\tiK}$.
 Note that each element $w\in{\sW}_{\tiJ}$ can be uniquely written as
 $w=w^{\tiK}w_{\tiK}$ with $w^{\tiK}\in{\sW}{}^{\tiK}_{\tiJ}$
 and $w_{\tiK}\in{\sW}_{\tiK}$, where $w^{\tiK}$ is a unique element
 of the smallest length in the right coset $w{\sW}_{\tiK}$.
 Moreover, we have
\[ {\sW}{}^{\tiK}_{\tiJ}= %
 \{\, w\in{\sW}_{\tiJ} \;|\; w(\al_k)>0 %
 \;\mbox{ for all }\; k\in{\bd K} \,\}, \]
 and
 ${\sW}{}^{\tiK}_{\tiJ}{\sW}{}^{\tiL}_{\tiK}={\sW}{}^{\tiL}_{\tiJ}$
 if ${\bd L}\subset{\bd K}\subset{\bd J}$.
 In addition, we set
\[ {\sDe}{}^{\tiK}_{\tiJ}:= %
 {\sDe}_{\tiJ}\setminus{\sDe}_{\tiK},\qquad %
 {\sDe}{}^{\tiK}_{\tiJ\pm}:= %
 {\sDe}{}^{\tiK}_{\tiJ}\cap{\sDe}_\pm. \]
 Note that ${\sDe}{}^{\emptyset}_{\tiJ}={\sDe}_{\tiJ}$
 and ${\sDe}{}^{\tiJ}_{\tiJ}=\emptyset$.

\begin{defn} 
 For each $w\in{\sW}_{\tiJ}$ and ${\bd K}\subset{\bd J}$,
 we set
\[ \De^{\tiK}_{\tiJ}\ft{(w,\pm)}:= %
 \la{w{\sDe}{}^{\tiK}_{\tiJ\pm}}\ra. \]
 We denote it simply by $\De_{\tiJ}\ft{(w,\pm)}$
 if ${\bd K}=\emptyset$, and by $\De^{\tiK}\ft{(w,\pm)}$
 if ${\bd J}={\SI}$.
\end{defn}

\begin{lem} 
 {\em(1)} The set $\De^{\tiK}_{\tiJ}\ft{(w,\pm)}$ is infinite
 if and only if ${\bd K}\subsetneq{\bd J}$.

 {\em(2)} The following formulas hold{\em:}
\begin{align}
 \De^{\tiK}_{\tiJ}\ft{(w,\pm)} &= 
 \De^{\tiK}_{\tiJ}\ft{(w^{\tiK},\pm)}, \\
 \De^{\tiK}_{\tiJ}\ft{(w^{\tiK},-)} &= 
 \vPhi(w^{\tiK})\amalg w^{\tiK}\De^{\tiK}_{\tiJ}\ft{(1,-)}, \\
 \re_{\tiJ+} &= 
 \De^{\tiK}_{\tiJ}\ft{(w,-)}\amalg w^{\tiK}\re_{\tiK+} %
 \amalg\De^{\tiK}_{\tiJ}\ft{(w,+)}.
\end{align}

 {\em(3)} If $\be\in w^{\tiK}\De_{\tiK+}$ and
 $\ga\in\De^{\tiK}_{\tiJ}\ft{(w,\pm)}$ satisfy
 $\be+\ga\in\De_+$, then $\be+\ga\in\De^{\tiK}_{\tiJ}\ft{(w,\pm)}$.

 {\em(4)} If a subset $B\subset\De_{\tiK+}$
 is biconvex in $\De_{\tiK+}$, then
 $\De^{\tiK}_{\tiJ}\ft{(w,\pm)}\amalg w^{\tiK}B$
 is biconvex in $\De_{\tiJ+}$.
\end{lem}

\begin{defn} 
 For each non-empty subset ${\bd J}\subset{\SI}$, we set
\begin{align*}
 \wt{\bol{\cl P}}_{\tiJ} &:= %
 \{\,(\scK,u,y) \;|\; \mbox{${\bd J}\supset{\bd K}$},\; %
 u\in{\sW}{}^{\tiK}_{\tiJ},\; y\in W_{\tiK}\,\}, \\
 \bol{\cl P}_{\tiJ} &:= %
 \{\,(\scK,u,y)\in\wt{\bol{\cl P}}_{\tiJ} \;|\; %
 \mbox{${\bd J}\supsetneq{\bd K}$}\,\}.
\end{align*}
 For each $(\scK,u,y)\in\wt{\bol{\cl P}}_{\tiJ}$,
 we define a subset
 $\nab_{\tiJ}\ft{(\tiK,u,y)}\subset\re_{\tiJ+}$ by setting
\[ \nab_{\tiJ}\ft{(\tiK,u,y)}:= %
 \De_{\tiJ}^{\tiK}\ft{(u,-)}\amalg u\vPhi_{\tiK}(y). \]
 Note that $\nab_{\tiJ}\ft{(\tiK,u,y)}=\vPhi_{\tiJ}(y)$
 if and only if ${\bd K}={\bd J}$ and that
 $\nab_{\tiJ}\ft{(\tiK,u,y)}=\De_{\tiJ}\ft{(u,-)}$
 if and only if ${\bd K}=\emptyset$.
\end{defn}

\begin{lem} 
 Let $(\scK_1,u_1,y_1)$ and $(\scK_2,u_2,y_2)$
 be elements of $\wt{\bol{\cl P}}_{\tiJ}$.
 Then the following two conditions are equivalent{\em:}
\[ {\mathrm(i)}\; \nab_{\tiJ}\ft{(\tiK_1,u_1,y_1)} %
 \;\dot{\subset}\;\nab_{\tiJ}\ft{(\tiK_2,u_2,y_2)}, %
 \qquad {\mathrm(ii)}\; %
 \mbox{${\bd K}_1\supset{\bd K}_2$},\;u_1\in u_2{\sW}_{\tiK_1}. \]
 In particular,
 $\nab_{\tiJ}\ft{(\tiK_1,u_1,y_1)} %
 \;\dot{\subset}\;\De_{\tiJ}\ft{(w,-)}$ for some
 $w\in{\sW}_{\tiJ}$ if and only if $u_1=w^{\tiK}$.
\end{lem}

\begin{thm} 
 {\em(1)} The assignment
 $(\scK,u,y)\,\mapsto\,\nab_{\tiJ}\ft{(\tiK,u,y)}$
 defines a bijective mapping between $\wt{\bol{\cl P}}_{\tiJ}$
 and ${\fr B}_{\tiJ}\amalg{\fr B}_{\tiJ}^\infty$, which maps
 $\bol{\cl P}_{\tiJ}$ onto ${\fr B}_{\tiJ}^\infty$.

 {\em(2)} The assignment $w\,\mapsto\,\De_{\tiJ}\ft{(w,-)}$
 defines a bijective mapping between ${\sW}_{\tiJ}$ and
 the set ${\fr M}_{\tiJ}$ of all maximal subsets of $\re_{\tiJ+}$
 {\em(}relative to the inclusion relation\/{\em)}
 which are convex in $\De_{\tiJ+}$.
 Moreover, ${\fr M}_{\tiJ}$ coincides with
 the set of all maximal subsets of $\re_{\tiJ+}$
 which are biconvex in $\De_{\tiJ+}$.
\end{thm}

\begin{prop}[\cite{jB}] 
 Let ${\bd K}$ be a proper subset of ${\bd J}$.
 Suppose that an element $\lam\in M_{\tiJ}$
 {\em(}see {\em Lemma 4.5)} satisfies $(\al_j\,|\,\lam)>0$
 for all $j\in{\bd J}\setminus{\bd K}$ and
 $(\al_k\,|\,\lam)=0$ for all $k\in{\bd K}$.
 Choose a reduced expression $t_{\lam}=\bs\ft{(1)}\cdots\bs\ft{(n)}$
 with $\bs\ft{(1)},\dots,\bs\ft{(n)}\in S_{\tiJ}$, and put
 $\bs\ft{(p)}=\bs\ft{(\ol{p})}$ for each $p\in{\bd N}$,
 where $\ol{p}$ is the unique positive integer such that
 $1\leq\ol{p}\leq n$ and $\ol{p}\equiv p\mod n$.
 Then we have $\bs\in{\cl W}_{\tiJ}^\infty$ and
 $\vPhi_{\tiJ}^\infty([\bs])=\De^{\tiK}_{\tiJ}\ft{(1,-)}$.
\end{prop}

\noindent
{\it Remark\/.}
 In case ${\bd J}={\SI}$ and ${\bd K}=\emptyset$,
 the above proposition was proved in \cite{jB}.

\begin{defn} 
 Let ${\bd K}$ be a proper subset of ${\bd J}$.
 We denote by $\ft{Z}^{\tiK}_{\tiJ}$
 the unique element of $W_{\tiJ}^\infty$ satisfying
 $\vPhi_{\tiJ}^\infty(\ft{Z}^{\tiK}_{\tiJ})= %
 \De^{\tiK}_{\tiJ}\ft{(1,-)}$ (see Proposition 5.8).
 Define a mapping
 $\chi_{\tiJ}\colon\bol{\cl P}_{\tiJ}\to W_{\tiJ}^\infty$
 by setting
 $\chi_{\tiJ}(\ft{(\tiK,u,y)})=uy.\ft{Z}^{\tiK}_{\tiJ}$
 for each $(\scK,u,y)\in\bol{\cl P}_{\tiJ}$.
\end{defn}

\begin{thm} 
 {\em(1)} For each element $x\in W_{\tiJ}$ and
 each proper subset ${\bd K}\subsetneq{\bd J}$,
 we have the following equality{\em:}
\begin{equation} 
 \vPhi_{\tiJ}^\infty(x.\ft{Z}^{\tiK}_{\tiJ})= %
 \nab_{\tiJ}\ft{(\tiK,\ol{x}^{\tiK},z_x)}
\end{equation}
 with a unique element $z_x\in W_{\tiK}$ such that
\begin{equation} 
 \vPhi_{\tiJ}((\ol{x}^{\tiK})^{-1}x)\cap\De_{\tiK+}= %
 \vPhi_{\tiK}(z_x).
\end{equation}

 {\em(2)} Let $\bol{\cl P}_{\tiJ}$ be the set defined
 in {\em Definition 5.5}, $\vPhi_{\tiJ}^\infty\colon %
 {W}_{\tiJ}^\infty\to{\fr B}_{\tiJ}^\infty$
 the injective mapping defined in {\em Definition 4.3},
 $\nab_{\tiJ}\colon\bol{\cl P}_{\tiJ}\to{\fr B}_{\tiJ}^\infty$
 the bijective mapping introduced in {\em Theorem 5.7}, and
 $\chi_{\tiJ}\colon\bol{\cl P}_{\tiJ}\to{W}_{\tiJ}^\infty$
 the mapping defined in {\em Definition 5.9}.
 Then both $\vPhi_{\tiJ}^\infty$ and $\chi_{\tiJ}$ are bijective
 and the following diagram is commutative{\em:}

\setlength{\unitlength}{0.5mm}
\begin{picture}(240,50)(-100,-5)
 \put(10,30){\makebox(20,10){${\fr B}_{\tiJ}^\infty$}}
 \put(-21,0){\makebox(20,10){${W}_{\tiJ}^\infty$}}
 \put(44,0){\makebox(20,10){$\bol{\cl P}_{\tiJ}$\;.}}
 \put(45,5){\vector(-1,0){48}}
 \put(-7,10){\vector(1,1){20}}
 \put(47,10){\vector(-1,1){20}}
 \put(17,0){$\chi_{\tiJ}$}
 \put(-8,20){$\vPhi_{\tiJ}^\infty$}
 \put(40,20){$\nab_{\tiJ}$}
%
\end{picture}

 {\em(3)} We have the following orbit decomposition{\em:}
\[ W_{\tiJ}^\infty= %
 \coprod_{{\bd K}\subsetneq{\bd J}}W_{\tiJ}.\ft{Z}^{\tiK}_{\tiJ}. \]
\end{thm}

\noindent
{\it Remark\/.}
 The existence and uniqueness of an element $z_x\in W_{\tiK}$
 satisfying (5.5) are guaranteed by Corollary 4.4 and
 the fact: for subsets $A$, $B$ and $C$ of $\De_+$ satisfying
 $B,\,C\subset A$, if $B$ is biconvex in $A$ then
 $B\cap C$ is biconvex in $C$.

%% file: sec6.tex
\section{Chains of infinite real biconvex sets}
 In this section, we consider chains
 $C_1\subset\cdots\subset C_n$ consisting of
 elements of ${\fr B}_{\tiJ}^\infty$, where
 ${\bd J}$ is a non-empty subset of ${\SI}$.
 Such chains are the key to
 the classification of convex orders on $\De_+$.

\begin{defn} 
 For each $w\in{\sW}_{\tiJ}$, we define
 $\wt{\bol{\cl P}}_{\tiJ}^w\subset\wt{\bol{\cl P}}_{\tiJ}$ and
 $\bol{\cl P}_{\tiJ}^w\subset\bol{\cl P}_{\tiJ}$ by setting
\begin{align*}
 \wt{\bol{\cl P}}_{\tiJ}^w &:= %
 \{\,(\scK,w^{\tiK},y)\in\wt{\bol{\cl P}}_{\tiJ}\;|\; %
 \vPhi_{\tiK}(y)\subset\De_{\tiK}\ft{(w_{\tiK},-)}\,\}, \\
 \bol{\cl P}_{\tiJ}^w &:= %
 \{\,(\scK,w^{\tiK},y)\in\wt{\bol{\cl P}}_{\tiJ}^w\;|\; %
 \mbox{${\bd J}\supsetneq{\bd K}$}\,\}.
\end{align*}
 Further we define $\bol{\cl Q}_{\tiJ}^w$ to be the set of
 all quadruplets $(\scK,\scL,y,z)$ consisting of subsets
 ${\bd L}\subset{\bd K}\subset{\bd J}$, $y\in W_{\tiK}$,
 and $z\in W_{\tiL}$ which satisfy the following conditions:
\[ {\mathrm Q(i)}\;\; %
 (\scK,w^{\tiK},y),\,(\scL,w^{\tiL},z)\in %
 \wt{\bol{\cl P}}_{\tiJ}^w;\qquad
 {\mathrm Q(ii)}\;\; %
 w^{\tiK}\vPhi_{\tiK}(y)\cap w^{\tiL}\De_{\tiL+} %
 \subset w^{\tiL}\vPhi_{\tiL}(z). \]
\end{defn}

\begin{prop} 
 Let $w$ be an element of ${\sW}_{\tiJ}$, and $(\scK,\scL,y,z)$
 a quadruplet consisting of subsets
 ${\bd L}\subset{\bd K}\subset{\bd J}$, $y\in W_{\tiK}$,
 and $z\in W_{\tiL}$.
 Then the following two conditions are equivalent{\em:}
\[ {\mathrm(i)}\;\; (\scK,\scL,y,z)\in\bol{\cl Q}_{\tiJ}^w;\qquad %
 {\mathrm(ii)}\;\; \nab_{\tiJ}\ft{(\tiK,w^{\tiK},y)}\subset %
 \nab_{\tiJ}\ft{(\tiL,w^{\tiL},z)}\subset\De_{\tiJ}\ft{(w,-)}. \]
 In particular, we have
 $\nab_{\tiJ}\ft{(\tiK,w^{\tiK},y)}\subset\De_{\tiJ}\ft{(w,-)}$
 if and only if $(\scK,w^{\tiK},y)\in\wt{\bol{\cl P}}_{\tiJ}^w$.
\end{prop}

\begin{pf}
 For each subset ${\bd S}\subset{\bd J}$, we have
\begin{equation} 
 \De_{\tiJ}\ft{(w,-)}=\De_{\tiJ}^{\tiS}\ft{(w^{\tiS},-)}\amalg %
 w^{\tiS}\De_{\tiS}\ft{(w_{\tiS},-)}
\end{equation}
 by (5.1) and (5.2). Moreover, by (5.1) we have
\begin{equation} 
 \De_{\tiJ}^{\tiK}\ft{(w^{\tiK},-)}\subset %
 \De_{\tiJ}^{\tiL}\ft{(w^{\tiL},-)}.
\end{equation}

 (i)$\Rightarrow$(ii)
 Suppose that $(\scK,\scL,y,z)\in\bol{\cl Q}_{\tiJ}^w$.
 Then both $\nab_{\tiJ}\ft{(\tiK,w^{\tiK},y)}$
 and $\nab_{\tiJ}\ft{(\tiL,w^{\tiL},z)}$ are subsets of
 $\De_{\tiJ}\ft{(w,-)}$ by Q(i) and (6.1).
 Thus we get
 $\nab_{\tiJ}\ft{(\tiK,w^{\tiK},y)}\subset %
 \nab_{\tiJ}\ft{(\tiL,w^{\tiL},z)}$ by (6.1), (6.2) and Q(ii).

 (ii)$\Rightarrow$(i)
 Suppose that
 $\nab_{\tiJ}\ft{(\tiK,w^{\tiK},y)}\subset %
 \nab_{\tiJ}\ft{(\tiL,w^{\tiL},z)}\subset\De_{\tiJ}\ft{(w,-)}$.
 By (6.1), both $(\scK,w^{\tiK},y)$ and $(\scL,w^{\tiL},z)$
 are elements of $\wt{\bol{\cl P}}_{\tiJ}^w$.
 Thus
 $w^{\tiK}\vPhi_{\tiK}(y)\cap w^{\tiL}\De_{\tiL+} %
 \subset w^{\tiL}\vPhi_{\tiL}(z)$ by (6.1) and (6.2).
\end{pf}

\begin{lem} 
 {\em(1)} Let $y$ be an element of $W_{\tiJ}$ such that
 $\vPhi_{\tiJ}(y)\subset\De_{\tiJ}\ft{(v,-)}$ for some
 $v\in{\sW}_{\tiJ}$. Then the following equality holds{\em:}
\[ \De_{\tiJ}\ft{(v,-)}= %
 \vPhi_{\tiJ}(y)\amalg y\De_{\tiJ}\ft{(\ol{y}^{-1}v,-)}. \]

 {\em(2)} Let ${\bd K}$ be a subset of ${\bd J}$.
 Then, for each $u\in{\sW}{}^{\tiK}_{\tiJ}$ and
 $(\scL,u',y)\in\wt{\bol{\cl P}}_{\tiK}$,
 the following equality holds{\em:}
\[ \De^{\tiK}_{\tiJ}\ft{(u,-)}\amalg %
 u\nab_{\tiK}\ft{(\tiL,u',y)}=\nab_{\tiJ}\ft{(\tiL,uu',y)}. \]
\end{lem}

\begin{pf}
 (1)\; Put
 $B:=y^{-1}\{\De_{\tiJ}\ft{(v,-)}\setminus\vPhi_{\tiJ}(y)\}$.
 Then $B$ is biconvex in $\De_{\tiJ+}$
 by Lemma 2.5 and Proposition 4.2.
 Since $\sharp{\vPhi_{\tiJ}(y)}<\infty$ we have
 $B\,\dot{=}\,\De_{\tiJ}\ft{(\ol{y}^{-1}v,-)}$
 by (iii) of Lemma 5.2(1).
 Thus we get $\De_{\tiJ}\ft{(\ol{y}^{-1}v,-)}\subset{B}$
 by Lemma 5.2(2), and hence $\De_{\tiJ}\ft{(\ol{y}^{-1}v,-)}=B$
 by Theorem 5.7(2).

 (2)\; Since $u\in{\sW}{}^{\tiK}_{\tiJ}$ and
 $\vPhi(u')\subset{\sDe}_{\tiK+}$, we have
 $u\vPhi(u')\subset{\sDe}_{\tiJ+}$. Hence we have
 $\vPhi(u)\amalg u\vPhi(u')=\vPhi(uu')$.
 Moreover, since $u'\in{\sW}_{\tiK}$ we have
 $u'\De^{\tiK}_{\tiJ}\ft{(1,-)}=\De^{\tiK}_{\tiJ}\ft{(1,-)}$.
 Therefore, by (5.2) we get
\begin{align*}
 (\mbox{the left hand side}) &= %
 \{\vPhi(u)\amalg u\De^{\tiK}_{\tiJ}\ft{(1,-)}\} %
 \amalg u\{u'\vPhi_{\tiL}(y)\amalg\vPhi(u')\amalg %
 u'\De^{\tiL}_{\tiK}\ft{(1,-)}\} \\
 &= \vPhi(uu')\amalg uu'\De^{\tiK}_{\tiJ}\ft{(1,-)}\amalg %
 uu'\vPhi_{\tiL}(y)\amalg uu'\De^{\tiL}_{\tiK}\ft{(1,-)} \\
 &= \vPhi(uu')\amalg uu'\De^{\tiL}_{\tiJ}\ft{(1,-)} %
 \amalg uu'\vPhi_{\tiL}(y) \\
 &= \De^{\tiL}_{\tiJ}\ft{(uu',-)}\amalg uu'\vPhi_{\tiL}(y) %
 = (\mbox{the right hand side}). \qed
\end{align*}
\renewcommand{\qed}{}
\end{pf}

\begin{prop} 
 {\em(1)} Let $(\scK,w^{\tiK},y)$ be an arbitrary element of
 $\wt{\bol{\cl P}}_{\tiJ}^w$, and put
 $C_1=\nab_{\tiJ}\ft{(\tiK,w^{\tiK},y)}$ and
 $v=\ol{y}^{-1}w_{\tiK}$.
 Then the following assignment
\[ B\,\longmapsto\,C=C_1\amalg w^{\tiK}yB \]
 defines a bijective mapping between the set
 $\{\,B\in{\fr B}_{\tiK}\amalg{\fr B}_{\tiK}^\infty \;|\; %
 B\subset\De_{\tiK}\ft{(v,-)}\,\}$ and the set
 $\{\,C\in{\fr B}_{\tiJ}\amalg{\fr B}_{\tiJ}^\infty \;|\; %
 C_1\subset C\subset\De_{\tiJ}\ft{(w,-)}\,\}$,
 which keeps the inclusion relation.

 {\em(2)} Under the condition in {\em(1)}, suppose that
 $B=\nab_{\tiK}\ft{(\tiL,v^{\tiL},g)}$ for some
 $(\scL,v^{\tiL},g)\in\wt{\bol{\cl P}}_{\tiK}^v$,
 and put $x=yv^{\tiL}g\in W_{\tiK}$.
 Then $C=\nab_{\tiJ}\ft{(\tiL,w^{\tiL},z_x)}$ with
 a unique element $z_x\in W_{\tiL}$ such that
\begin{equation} 
 \vPhi_{\tiK}((\ol{x}^{\tiL})^{-1}x)\cap\De_{\tiL+}= %
 \vPhi_{\tiL}(z_x).
\end{equation}
\end{prop}

\begin{pf}
 (1)\; By the definition of $\wt{\bol{\cl P}}_{\tiJ}^w$
 and Lemma 6.3(1), we have
\begin{equation} 
 \De_{\tiK}\ft{(w_{\tiK},-)}= %
 \vPhi_{\tiK}(y)\amalg y\De_{\tiK}\ft{(v,-)}.
\end{equation}
 Thus, by (6.1) and (6.4) we get
\begin{equation} 
 \De_{\tiJ}\ft{(w,-)}=C_1\amalg w^{\tiK}y\De_{\tiK}\ft{(v,-)}.
\end{equation}

 Let $B$ be an element of
 ${\fr B}_{\tiK}\amalg{\fr B}_{\tiK}^\infty$
 such that $B\subset\De_{\tiK}\ft{(v,-)}$,
 and put $C:=C_1\amalg w^{\tiK}yB$.
 Then $C_1\subset C\subset\De_{\tiJ}\ft{(w,-)}$ by (6.5).
 Put $D:=\vPhi_{\tiK}(y)\amalg yB$.
 Then we see that
 $C=\De^{\tiK}_{\tiJ}\ft{(w^{\tiK},-)}\amalg w^{\tiK}D$ and
 that $D$ is biconvex in $\De_{\tiK+}$ by Lemma 2.5.
 Thus $C$ is biconvex in $\De_{\tiJ+}$ by Lemma 5.4(4),
 and hence the mapping is well-defined.

 The injectivity is clear. Let us check the surjectivity.
 Suppose that $C\in{\fr B}_{\tiJ}\amalg{\fr B}_{\tiJ}^\infty$
 satisfies $C_1\subset C\subset\De_{\tiJ}\ft{(w,-)}$.
 By (6.5), there exists a unique subset
 $B\subset\De_{\tiK}\ft{(v,-)}$ such that $C=C_1\amalg w^{\tiK}yB$.
 Then we see that
 $C\cap w^{\tiK}\De_{\tiK+}=w^{\tiK}(\vPhi_{\tiJ}(y)\amalg yB)$
 and this set is biconvex in $w^{\tiK}\De_{\tiK+}$.
 Thus $B\in{\fr B}_{\tiK}\amalg{\fr B}_{\tiK}^\infty$ by Lemma 2.5.
 Therefore the mapping is surjective.

 (2)\; In case ${\bd K}={\bd L}$, we have
 $B=\vPhi_{\tiL}(g)$ and  $z_x=x=yg$.
 Since $yB\subset\De_{\tiL+}$ we have
 $\vPhi_{\tiL}(y)\amalg yB=\vPhi_{\tiL}(yg)$.
 Thus we get
\begin{align*}
 C &= \De^{\tiL}_{\tiJ}\ft{(w^{\tiL},-)}\amalg %
 w^{\tiL}\vPhi_{\tiL}(y)\amalg w^{\tiL}yB \\
 &= \De^{\tiL}_{\tiJ}\ft{(w^{\tiL},-)}\amalg %
 w^{\tiL}\vPhi_{\tiL}(yg)
 = \nab_{\tiJ}\ft{(\tiL,w^{\tiL},z_x)}.
\end{align*}

 In case ${\bd K}\supsetneq{\bd L}$,
 since $yB\subset\De_{\tiK+}$ we have
 $\vPhi_{\tiK}(y)\amalg yB= %
 \vPhi_{\tiK}^{\infty}(y.v^{\tiL}g.\ft{Z}^{\tiL}_{\tiK})$
 by Proposition 2.12 and Theorem 5.10(2).
 By Theorem 5.10(1) we have
\[ \vPhi_{\tiK}(y)\amalg yB= %
 \vPhi_{\tiK}^{\infty}(x.\ft{Z}^{\tiL}_{\tiK})= %
 \nab_{\tiK}\ft{(\tiL,\ol{x}^{\tiL},z_x)}. \]
 Thus, by Lemma 6.3(2) we get
\begin{align*}
 C &= \De_{\tiJ}^{\tiK}\ft{(w^{\tiK},-)}\amalg %
 w^{\tiK}\vPhi_{\tiK}(y)\amalg w^{\tiK}yB \\
 &= \De_{\tiJ}^{\tiK}\ft{(w^{\tiK},-)}\amalg %
 w^{\tiK}\nab_{\tiK}\ft{(\tiL,\ol{x}^{\tiL},z_x)}
 = \nab_{\tiJ}\ft{(\tiL,w^{\tiK}\ol{x}^{\tiL},z_x)}.
\end{align*}
 Here, by Lemma 5.6 we have
 $w^{\tiK}\ol{x}^{\tiL}=w^{\tiL}$ since
 $C\subset\De_{\tiJ}\ft{(w,-)}$.
\end{pf}

\begin{cor} 
 Choose an element $(\scK,w^{\tiK},y)\in\wt{\bol{\cl P}}_{\tiJ}^w$,
 and put $v=\ol{y}^{-1}w_{\tiK}$. Further, choose an element
 $(\scL,v^{\tiL},g)\in\wt{\bol{\cl P}}_{\tiK}^v$,
 and put $x=yv^{\tiL}g\in W_{\tiK}$.
 Then $(\scK,\scL,y,z_x)$ is an element of $\bol{\cl Q}_{\tiJ}^w$,
 where $z_x$ is the unique element of $W_{\tiL}$ satisfying
 the equality {\em(6.3)}.
 Moreover, each element of $\bol{\cl Q}_{\tiJ}^w$
 is constructed by applying the above procedure suitably.
\end{cor}

\begin{pf}
 This follows immediately from Proposition 6.2 and Proposition 6.4.
\end{pf}

\begin{defn} 
 (1) For each positive integer $n\in{\bd N}$, we define
 $\bol{\cl C}_n\De_{\tiJ}\ft{(w,-)}$ to be the set of all
 $(C_0,\,C_1,\dots,C_n)$ consisting of $C_0=\emptyset$ and
 $C_1,\dots,C_n\in{\fr B}_{\tiJ}^\infty$ such that
 $\emptyset=C_0\subset %
 C_1\subset\cdots\subset C_n=\De_{\tiJ}\ft{(w,-)}$ and
 $\sharp(C_{i'}\setminus C_i)=\infty$ if $i<i'$.
 We denote an element $(C_0,\,C_1,\dots,C_n)$ of
 $\bol{\cl C}_n\De_{\tiJ}\ft{(w,-)}$ by $C_\bullet$.
 Further we put
\[\bol{\cl C}\De_{\tiJ}\ft{(w,-)}:=\coprod_{n=1}^\infty
 \bol{\cl C}_n\De_{\tiJ}\ft{(w,-)}. \]

 (2) For each positive integer $n\in{\bd N}$, we set
\[ \bol{\cl C}_n{\bd J}:=\{\, {\bd K}_\bul= %
 ({\bd K}_0,{\bd K}_1,\dots,{\bd K}_n) \;|\; %
 {\bd J}={\bd K}_0\supsetneq{\bd K}_1 %
 \supsetneq\cdots\supsetneq{\bd K}_n=\emptyset \,\}. \]
 Note that $\bol{\cl C}_n{\bd J}=\emptyset$ if $n>\sharp{\bd J}$.
 Put
\[ \bol{\cl C}{\bd J}:= %
 \coprod_{n=1}^{\sharp{\bd J}}\bol{\cl C}_n{\bd J}. \]
 For each $n=1,\dots,\sharp{\bd J}$ and
 ${\scK}_\bullet\in\bol{\cl C}_n{\bd J}$, we put
\begin{align*}
 W_{\tiK_\bul} &:= %
 W_{\tiK_1} \times\cdots\times W_{\tiK_n}, \\
 {\cl W}_{\tiK_\bul}^\infty &:= {\cl W}_{\tiK_0}^\infty %
 \times\cdots\times {\cl W}_{\tiK_{n-1}}^\infty, \\
 {\fr B}_{\tiK_\bul}^\infty &:= {\fr B}_{\tiK_0}^\infty %
 \times\cdots\times{\fr B}_{\tiK_{n-1}}^\infty.
\end{align*}
 Denote an element
 $(y_1,\dots,y_n)\in W_{\tiJ_\bullet}$ by $y_\bullet$,
 an element $({\bs}_0,\dots,{\bs}_{n-1})\in{\cl W}_{\tiJ_\bul}^\infty$
 by ${\bs}_\bul$, and an element
 $(B_0,\dots,B_{n-1})\in{\fr B}_{\tiJ_\bul}^\infty$
 by $B_\bul$.
\end{defn}

\begin{prop} 
 Let $C_\bul$ be an arbitrary element of
 $\bol{\cl C}_n\De_{\tiJ}\ft{(w,-)}$.
 Then there exist a unique ${\scK}_\bul\in\bol{\cl C}_n{\bd J}$
 and a unique $y_\bul\in W_{\tiK_\bullet}$ such that
 $\nab_{\tiJ}\ft{(\tiK_i,w^{\tiK_i},y_i)}=C_i$ for $i=1,\dots,n$.
 In particular, we have $n\leq\sharp{\bd J}$, and hence
 $\bol{\cl C}_n\De_{\tiJ}\ft{(w,-)}=\emptyset$ if $n>\sharp{\bd J}$.
\end{prop}

\begin{pf}
 For each $i=1,\dots,n$, there exists a unique element
 $(\scK_i,w^{\tiK_i},y_i)\in\bol{\cl P}_{\tiJ}$ such that
 $C_i=\nab_{\tiJ}\ft{(\tiK_i,w^{\tiK_i},y_i)}$
 by Lemma 5.6 and Theorem 5.7(1).
 Since $C_i\subset C_{i'}$ and
 $\sharp(C_{i'}\setminus C_i)=\infty$ for each $i<i'$,
 we have ${\bd K}_i\supsetneq{\bd K}_{i'}$ by Lemma 5.6.
 In addition, we have ${\bd K}_n=\emptyset$ since
 $C_n=\De_{\tiJ}\ft{(w,-)}$.
\end{pf}

\begin{defn} 
 For each $n=1,\dots,\sharp{\bd J}$, we define
 $\bol{\cl S}_n\ft{(\scJ,w)}$ to be the set of
 all pairs $(\scK_\bul,y_\bul)$ consisting of
 $\scK_\bul\in\bol{\cl C}_n{\bd J}$ and
 $y_\bul\in W_{\tiK_\bul}$ such that
 $(\scK_{i-1},\scK_i,y_{i-1},y_i)\in\bol{\cl Q}_{\tiJ}^w$
 for each $i=1,\dots,n$, where $y_0=1$. Further we put
\[ \bol{\cl S}\ft{(\scJ,w)}:= %
 \coprod_{n=1}^{\sharp{\bd J}}\bol{\cl S}_n\ft{(\scJ,w)}. \]

 For each $(\scK_\bul,y_\bul)\in\bol{\cl S}_n\ft{(\scJ,w)}$
 and $i=1,\dots,n$, we set
\[ C_i\ft{(\tiK_\bul,y_\bul)}:= %
 \nab_{\tiJ}\ft{(\tiK_i,w^{\tiK_i},y_i)}, \]
 and put
\[ C_\bul\ft{(\tiK_\bul,y_\bul)}:= %
 (C_0\ft{(\tiK_\bul,y_\bul)},C_1\ft{(\tiK_\bul,y_\bul)}, %
 \dots,C_n\ft{(\tiK_\bul,y_\bul)}), \]
 where $C_0\ft{(\tiK_\bul,y_\bul)}:=\emptyset$.
\end{defn}

\noindent
{\it Remark\/.}
 An element of $\bol{\cl S}\ft{(\scJ,w)}$ is constructed by
 applying Corollary 6.5 recursively.

\begin{thm} 
 For each non-empty subset ${\bd J}\subset{\SI}$ and
 each element $w\in{\sW}_{\tiJ}$, the assignment
 $(\scK_\bul,y_\bul)\mapsto C_\bul\ft{(\tiK_\bul,y_\bul)}$
 defines a bijective mapping between
 $\bol{\cl S}\ft{(\scJ,w)}$ and
 $\bol{\cl C}\De_{\tiJ}\ft{(w,-)}$,
 which maps $\bol{\cl S}_n\ft{(\scJ,w)}$ onto
 $\bol{\cl C}_n\De_{\tiJ}\ft{(w,-)}$
 for each $n=1,\dots,\sharp{\bd J}$.
\end{thm}

\begin{pf}
 Let $(\scK_\bul,y_\bul)$ be an element of
 $\bol{\cl S}_n\ft{(\scJ,w)}$.
 Since
 $(\scK_{i-1},\scK_i,y_{i-1},y_i)\in\bol{\cl Q}_{\tiJ}^w$
 for each $i=1,\dots,n$, we have
 $C_{i-1}\ft{(\tiK_\bul,y_\bul)}\subset %
 C_i\ft{(\tiK_\bul,y_\bul)}\subset\De_{\tiJ}\ft{(w,-)}$
 by Proposition 6.2.
 Since ${\bd K}_n=\emptyset$ we have
 $C_n\ft{(\tiK_\bul,y_\bul)}=\De_{\tiJ}\ft{(w,-)}$.
 Moreover, since
 ${\bd K}_{i-1}\supsetneq{\bd K}_i$ for each $i=1,\dots,n$,
 we have $\sharp(C_i\ft{(\tiK_\bul,y_\bul)}\setminus %
 C_{i-1}\ft{(\tiK_\bul,y_\bul)})=\infty$ by Lemma 5.6.
 Thus $C_\bul\ft{(\tiK_\bul,y_\bul)}$ is an element of
 $\bol{\cl C}_n\De_{\tiJ}\ft{(w,-)}$.
 The injectivity follows from Theorem 5.7(1), and
 the surjectivity follows from Proposition 6.7 and Proposition 6.2.
\end{pf}

\begin{prop} 
 For each $n=1,\dots,\sharp{\bd J}$ and
 $(\scK_\bul,y_\bul)\in\bol{\cl S}_n\ft{(\scJ,w)}$,
 there exists a unique
 $B_\bullet\ft{(\tiK_\bul,y_\bul)}= %
 (B_0\ft{(\tiK_\bul,y_\bul)},\dots,B_{n-1}\ft{(\tiK_\bul,y_\bul)})
 \in{\fr B}_{\tiK_\bullet}^\infty$ which satisfies
 the following equality{\em:}
\begin{equation} 
 C_i\ft{(\tiK_\bul,y_\bul)}= %
 C_{i-1}\ft{(\tiK_\bul,y_\bul)}\amalg %
 w^{\tiK_{i-1}}y_{i-1}B_{i-1}\ft{(\tiK_\bul,y_\bul)}
\end{equation}
 for $i=1,\dots,n$, where $y_0=1$.
\end{prop}

\begin{pf}
 By Theorem 6.9, we have
\[ C_{i-1}\ft{(\tiK_\bul,y_\bul)}\subset %
 C_i\ft{(\tiK_\bul,y_\bul)}\subset\De_{\tiJ}\ft{(w,-)}, %
 \quad\mbox{and}\quad \sharp(C_i\ft{(\tiK_\bul,y_\bul)} %
 \setminus C_{i-1}\ft{(\tiK_\bul,y_\bul)})=\infty \]
 for each $i=1,\dots,n$.
 Thus, by Proposition 6.4 there exists a unique
 $B_{i-1}\ft{(\tiK_\bul,y_\bul)}\in{\fr B}_{\tiK_{i-1}}^\infty$
 which satisfies the equality (6.6).
\end{pf}

\begin{defn} 
 For each $n=1,\dots,\sharp{\bd J}$, we define
 $\bol{\cl T}_n\ft{(\scJ,w)}$ to be the set of all triplets
 $(\scK_\bul,y_\bul,\bs_\bul)$ consisting of
 $(\scK_\bul,y_\bul)\in\bol{\cl S}_n\ft{(\scJ,w)}$ and
 $\bs_\bul\in{\cl W}_{\tiK_\bul}^\infty$ such that
\begin{equation} 
 B_{i-1}\ft{(\tiK_\bul,y_\bul)}= %
 \vPhi_{\tiK_{i-1}}^\infty([\bs_{i-1}])
\end{equation}
 for $i=1,\dots,n$. Further we put
\[ \bol{\cl T}\ft{(\scJ,w)}:= %
 \coprod_{n=1}^{\sharp{\bd J}}\bol{\cl T}_n\ft{(\scJ,w)}. \]
\end{defn}

\noindent
{\it Remark\/.}
 The existence of an element
 $\bs_{i-1}\in{\cl W}_{\tiK_{i-1}}^\infty$ satisfying (6.7)
 is guaranteed by Theorem 5.10(2).

%% file: sec7.tex
\section{The classification of convex orders}
 In this section, we classify all convex orders on $\De_+$ and
 give a general method of the construction of convex orders on
 $\De_+$ for an arbitrary untwisted affine Lie algebra.
 To solve the problems, we have to classify and to construct
 all convex orders on the subset $\De_{\tiJ+}\subset\De_+$
 for each non-empty subset ${\bd J}\subset{\SI}$.

\begin{lem} 
 Let $u$ be an element of ${\sW}{}^{\tiJ}$.

 {\em(1)} Let $\preq$ be a total order on a subset
 $B\subset u\De_{\tiJ+}$.
 Then the conditions $\mathrm{CO(ii)}_{\De_+}$ and
 $\mathrm{CO(ii)}_{u\De_{\tiJ+}}$ concerning $\preq$ are equivalent.

 {\em(2)} Let $\preq$ be a convex order on
 a subset $B\subset\De_+$, and put $C:=B\cap u\De_{\tiJ+}$.
 Then the restriction $\preq_{\ti{C}}$ is a convex order.
\end{lem}

\begin{pf}
 (1) Suppose that $\be\in B$ and $\ga\in\De_+\setminus B$
 satisfy $\be+\ga\in B$.
 Then $\ga\in u\De_{\tiJ+}\setminus B$ by Lemma 5.4(3).
 Thus $\mbox{CO(ii)}_{u\De_{\tiJ+}}$ implies
 $\mbox{CO(ii)}_{\De_+}$. The inverse is clear.

 (2) By (1), it suffices to check the conditions CO(i) and
 $\mathrm{CO(ii)}_{u\De_{\tiJ+}}$. The condition CO(i) is clear.
 Suppose that $\be\in C$ and $\ga\in u\De_{\tiJ+}\setminus C$
 satisfy $\be+\ga\in C$.
 Since $u\De_{\tiJ+}\setminus C=u\De_{\tiJ+}\setminus B$ we have
 $\be\pr\be+\ga$ by the property $\mbox{CO(ii)}_{\De_+}$ of $\preq$.
 Thus $\preq_{\ti{C}}$ satisfies $\mathrm{CO(ii)}_{u\De_{\tiJ+}}$.
\end{pf}

\noindent
{\it Remark\/.}
 Thanks to Lemma 7.1(1), when we check the property
 $\mbox{CO(ii)}_{\De_+}$ of a total order $\preq$
 on a subset $B\subset u\De_{\tiJ+}$, it is enough to check
 the condition $\mbox{CO(ii)}_{u\De_{\tiJ+}}$.

\begin{lem} 
 Let $\preq$ be a convex order on $\De_{\tiJ+}$.

 {\em(1)} If $\be\in\De_{\tiJ+}$ satisfies $m\de\pr \be\pr n\de$
 for some $m,\,n\in{\bd N}$, then $\be\in\pim$.

 {\em(2)} If $\vep\in{\sDe}_{\tiJ}$ satisfies
 $m\de+\vep\in\De_{\tiJ+}$ and $m\de+\vep\pr\de$
 {\em(resp. $\de\pr m\de+\vep$)} for some $m\in{\bd Z}_{\geq0}$,
 then $m\de+\vep\pr(m+1)\de+\vep\pr\de$
 {\em(resp. $\de\pr(m+1)\de+\vep\pr m\de+\vep$)}.

 {\em(3)} For each $\vep\in{\sDe}_{\tiJ}$, we have
 either $\la{\vep}\ra\pr\pim$ or $\pim\pr\la{\vep}\ra$.
\end{lem}

\begin{pf}
 (1) Suppose that
 $m\de\pr\be\pr n\de$ for some $\be\in\re_{\tiJ+}$.
 Then $m\de\pr mn\de+\be\pr\be\pr mn\de+\be\pr n\de$
 by the property CO(i). This is a contradiction.

 (2) Clear.

 (3) By (2), we have either
 $\la{\vep}\ra\pr\{\de\}$ or $\{\de\}\pr\la{\vep}\ra$.
 Hence, by (1) we get
 either $\la{\vep}\ra\pr\pim$ or $\pim\pr\la{\vep}\ra$.
\end{pf}

\begin{thm} 
 Let us choose an element $w\in{\sW}_{\tiJ}$,
 a convex order $\preq_-$ on $\De_{\tiJ}\ft{(w,-)}$,
 a total order $\preq_0$ on $\pim$, and
 an opposite convex order $\preq_+$ on $\De_{\tiJ}\ft{(w,+)}$.
 Define a total order $\preq$ on $\De_{\tiJ+}$ by extending
 $\preq_-$, $\preq_0$, and $\preq_+$ to
 $\De_{\tiJ+}=\De_{\tiJ}\ft{(w,-)}\amalg\pim\amalg\De_{\tiJ}\ft{(w,+)}$
 as follows{\em:}
\[ \De_{\tiJ}\ft{(w,-)}\,\pr\,\pim\,\pr\,\De_{\tiJ}\ft{(w,+)}. \]
 Then $\preq$ is a convex order on $\De_{\tiJ+}$.
 Moreover, each convex order on $\De_{\tiJ+}$ is constructed by
 applying the above procedure suitably.
\end{thm}

\begin{pf}
 Let $\preq$ be a total order on $\De_{\tiJ+}$
 constructed by applying the above procedure.
 Then it is easy to see that $\preq$ is a convex order on
 $\De_{\tiJ+}$ by the definition of convex orders on
 $\De\ft{(w,-)}$ and opposite convex orders on $\De\ft{(w,+)}$.

 Conversely, suppose that
 $\preq$ is a convex order on $\De_{\tiJ+}$.
 By Lemma 7.2(1), there exist unique subsets
 $B_-,\,B_+\subset\re_{\tiJ+}$ satisfying
 $B_-\amalg B_+=\re_{\tiJ+}$ and $B_-\pr\pim\pr B_+$.
 Let $\preq_\pm$ be the restriction of $\preq$ to $B_\pm$.
 Then $\preq_-$ is a convex order on $B_-$ and $\preq_+$
 is an opposite convex order on $B_+$ by Lemma 3.4(1)(2).
 Therefore, to complete the proof, it is enough to show that
 $B_\pm=\De_{\tiJ}\ft{(w,\pm)}$ with a unique $w\in{\sW}_{\tiJ}$.
 Put
 $P_\pm:=\{\,\vep\in{\sDe}_{\tiJ}\;|\;\la{\vep}\ra\subset B_\pm\,\}$.
 By Lemma 7.2(3), we have
\[ \mbox{(i)}\;\; B_\pm = \la{P_\pm}\ra, \qquad
 \mbox{(ii)}\;\; P_-\,\amalg\,P_+ = {\sDe}_{\tiJ}. \]
 Moreover, both subsets $B_-,\,B_+\subset\re_{\tiJ+}$
 are convex in $\De_{\tiJ+}$ by Lemma 3.4(1)(2).
 This fact implies that
\[ \mbox{(iii)}\;\; \vep,\,\eta\in P_\pm,\;\vep+\eta\in{\sDe}_{\tiJ} %
 \,\Lra\, \vep+\eta\in P_\pm, \qquad %
 \mbox{(iv)}\;\; P_\pm\cap(-P_\pm)=\emptyset. \]
 By (ii) and (iv), we have $-P_-=P_+$, and hence
\[ \mbox{(v)}\;\; {\sDe}_{\tiJ} = P_-\amalg(-P_-). \]
 Thanks to (iii) and (v), there exists a unique
 $w\in{\sW}_{\tiJ}$ such that $P_-=w{\sDe}_{\tiJ-}$
 by the theory of classical root systems (cf. \cite{nB}).
 Thus we get $B_\pm=\De_{\tiJ}\ft{(w,\pm)}$ by (i).
\end{pf}

\noindent
{\it Remark\/.}
 Thanks to Theorem 7.3, to complete the classification of
 convex orders on $\De_{\tiJ+}$, it suffices to classify all
 convex orders on $\De_{\tiJ}\ft{(w,-)}$ for each $w\in{\sW}_{\tiJ}$
 since $\De_{\tiJ}\ft{(w,+)}=\De_{\tiJ}\ft{(ww_\circ,-)}$,
 where $w_\circ$ is the longest element of ${\sW}_{\tiJ}$.

\begin{lem} 
 Let $\preq$ be a total order on a subset
 $B\subset\re_+$ with the property {\em CO(ii)}.

 {\em(1)} Suppose that
 $m\de+\vep$ and $(m+k)\de+\vep$ are elements of $B$ with
 $m\in{\bd Z}_{\geq0}$, $k\in{\bd N}$, and $\vep\in{\sDe}$.
 Then $m\de+\vep\pr(m+k)\de+\vep$.

 {\em(2)} The order $\preq$ is a well-order on $B$.

 {\em(3)} Suppose that $B=\la{P}\ra$ for some $P\subset{\sDe}$.
 Then there is not the maximum element of $B$ relative to $\preq$.
\end{lem}

\begin{pf}
 (1) Since $B\subset\pre$ we have $k\de\in\De_+\setminus B$.
 Thus we get $m\de+\vep\pr(m+k)\de+\vep$
 by the property CO(ii) of $\preq$.

 (2) Let $A$ be a non-empty subset of $B$.
 Define a subset $\ol{A}\subset{\sDe}$ by setting
 $\ol{A}=\{\,\ol{\al}\,|\,\al\in A\,\}$.
 Then $A=\coprod_{\vep\in\ol{A}}(\la{\vep}\ra\cap A)$.
 For each $\vep\in\ol{A}$, by (1) there exists the minimum element
 $\xi_\vep$ of $\la{\vep}\ra\cap A$ relative to $\preq$.
 Since $\ol{A}$ is finite, there exists the minimum element
 $\xi_\ast$ of $\{\,\xi_\vep\,|\,\vep\in\ol{A}\,\}$.
 Then $\xi_\ast$ is the minimum element of $A$.

 (3) Let $\xi$ be an element of $B$.
 By the definition of $\la{P}\ra$, we have
 $\xi\in\la{\vep}\ra$ for some unique $\vep\in P$.
 Then we have $\de+\xi\in\la{P}\ra$ and $\xi\pr\de+\xi$ by (1).
 Hence there is not the maximum element of $B$ relative to $\preq$.
\end{pf}

\begin{prop} 
 Let $w$ be an arbitrary element of ${\sW}_{\tiJ}$.
 Then, for each convex order $\preq$ on $\De_{\tiJ}\ft{(w,-)}$,
 there exists a unique
 $C_\bul\ft{(\preq)}\in\bol{\cl C}_n\De_{\tiJ}\ft{(w,-)}$
 with a unique positive integer $n\leq\sharp{\bd J}$ such that
\begin{itemize}
 \item[(i)] the restriction of $\preq$ to
 $R_i\ft{(\preq)}:= %
 C_i\ft{(\preq)}\setminus C_{i-1}\ft{(\preq)}$
 is of 1-row type for each $i=1,\dots,n${\em;}
 \item[(ii)] $R_i\ft{(\preq)}\pr R_{i'}\ft{(\preq)}$
 if and only if $i<i'$ for $i,\,i'=1,\dots,n$.
\end{itemize}
 In particular, we have
 $\De_{\tiJ}\ft{(w,-)}=\amalg_{i=1}^nR_i\ft{(\preq)}$,
 and hence $\preq$ is of n-row type.
\end{prop}

\begin{pf}
 By combining Lemma 3.2(2) with Lemma 7.4(2)(3), we see that
 $\De_{\tiJ}\ft{(w,-)}$ uniquely decomposes into
 the form $\amalg_{\lam\in\Lam}R_\lam\ft{(\preq)}$ with a uniquely
 determined (up to order isomorphism) well-ordered
 non-empty set $(\Lam,\wt{\preq})$ having the properties
 (i) and (ii) in Lemma 3.2(2).
 Put $C_\lam\ft{(\preq)}= %
 \amalg_{\mu\,\wt{\preq}\,\lam}R_\mu\ft{(\preq)}$
 for each $\lambda\in\Lambda$.
 Then $C_\lam\ft{(\preq)}$ is biconvex in $\De_{\tiJ+}$
 by Lemma 3.4(2). Since
 $\sharp(C_{\lam'}\ft{(\preq)}\setminus %
 C_\lam\ft{(\preq)})=\infty$ for each $\lam\,\wt{\pr}\,\lam'$,
 we get $\sharp\Lam\leq\sharp{\bd J}$ by Proposition 6.7.
 Hence we may identify $\Lam$ with the set $\{1,\cdots,n\}$, and
 $\wt{\preq}$ with the usual order $\leq$ on $\{1,\cdots,n\}$,
 where $n=\sharp\Lam$.
\end{pf}

\begin{defn} 
 We denote by $\bol{\cl{CO}}_n\De_{\tiJ}\ft{(w,-)}$ the set of
 all $n$-row type convex orders on $\De_{\tiJ}\ft{(w,-)}$, and
 by $\bol{\cl{CO}}\De_{\tiJ}\ft{(w,-)}$ the set of
 all convex orders on $\De_{\tiJ}\ft{(w,-)}$.
 Note that
\[ \bol{\cl{CO}}\De_{\tiJ}\ft{(w,-)}= %
 \coprod_{n=1}^{\sharp{\bd J}}\bol{\cl{CO}}_n\De_{\tiJ}\ft{(w,-)}. \]
\end{defn}

\begin{prop} 
 Let $C$ and $C_1$ be elements of
 ${\fr B}_{\tiJ}\amalg{\fr B}_{\tiJ}^\infty$
 satisfying $C=C_1\amalg uyB$ and
 $C_1=\nab_{\tiJ}\ft{({\tiK,u,y})}$ with
 $(\scK,u,y)\in\wt{\bol{\cl P}}_{\tiJ}$ and
 $B\in{\fr B}_{\tiK}\amalg{\fr B}_{\tiK}^\infty$.

 {\em(1)} If there exist a convex order $\preq_1$ on $C_1$
 and a convex order $\preq'$ on $B$, then there exists
 a unique convex order $\preq$ on $C$ such that
\begin{itemize}
 \item[(i)] $C_1\pr uyB$,
 \item[(ii)] the restriction of $\preq$ to $C_1$ is $\preq_1$,
 \item[(iii)] $uy(\xi)\preq uy(\eta)$ if and only if
 $\xi\preq'\eta$ for $\xi,\,\eta\in B$.
\end{itemize}

 {\em(2)} If there exists a convex order $\preq$ on $C$ with
 the above property {\em(i)}, then there exists a unique
 convex order $\preq'$ on $B$ with
 the above property {\em(iii)}.
\end{prop}

\begin{pf}
 Put $D:=C\cap u\De_{\tiK+}$. Then we have
\begin{align} 
 D &= u\vPhi_{\tiK}(y)\amalg uyB, \\
 C &= \De^{\tiK}_{\tiJ}\ft{(u,-)}\amalg D, \\
 \De_{\tiJ+}\setminus C &= %
 \{u\De_{\tiK+}\setminus D\}\amalg\De^{\tiK}_{\tiJ}\ft{(u,+)}.
\end{align}

 (1) Since the properties (i)--(iii) determine
 a unique total order $\preq$ on $C$, it suffices to
 check that $\preq$ is a convex order on $C$.
 We first show that
 the restriction $\preq_{\ti{D}}$ is a convex order.
 Indeed, the restriction of $\preq$ to $u\vPhi(y)$ is
 a convex order by Lemma 7.1(2), and hence $\preq_{\ti{D}}$
 is a convex order by (7.1) and Proposition 3.6(1).
 Let us check the property CO(i) of $\preq$.
 Suppose that $\be\in C$ and $\ga\in C$ satisfy
 $\be\pr\ga$ and $\be+\ga\in\De_{\tiJ+}$.
 By the properties (ii) and (iii), it is enough to
 consider the case: $\be\in C_1$ and $\ga\in uyB$.
 Moreover, we may assume that $\be\in\De^{\tiK}_{\tiJ}\ft{(u,-)}$
 since $\preq_{\ti{D}}$ is a convex order.
 Then $\be+\ga\in\De^{\tiK}_{\tiJ}\ft{(u,-)}$ by Lemma 5.4(3).
 Thus we get $\be\pr\be+\ga$ by the property CO(ii) of $\preq_1$,
 and hence $\be\pr\be+\ga\pr\ga$ by the property (i).
 We next check the property CO(ii) of $\preq$.
 Suppose that $\be\in C$ and
 $\ga\in\De_{\tiJ+}\setminus C$ satisfy $\be+\ga\in C$.
 By the property CO(ii) of $\preq_1$,
 we may assume that $\be\in uyB$.
 Then we have $\ga\in u\De_{\tiK+}\setminus D$
 by (7.3) and Lemma 5.4(3), and hence $\be+\ga\in D$.
 Thus $\be\pr\be+\ga$ since $\preq_{\ti{D}}$ is a convex order.

 (2) By Lemma 7.1(2), the restriction $\preq_{\ti{D}}$ is
 a convex order such that $u\vPhi_{\tiK}(y)\pr_{\ti{D}}uyB$.
 Thus this follows from (7.1) and Proposition 3.6(2).
\end{pf}

\begin{defn} 
 For each $n=1,\dots,\sharp{\bd J}$ and
 $(\scK_\bul,y_\bul,\bs_\bul)\in\bol{\cl T}_n\ft{(\scJ,w)}$
 (see Definition 6.11), we define a total order
 $\preq_\bul\!\!\ft{(\tiK_\bul,y_\bul,\bs_\bul)}=\preq_\bul$
 on $\De_{\tiJ}\ft{(w,-)}$ by applying
 the following procedure Step1--3.

 Step 1.
 For each $i=1,\dots,n$, put
\[ R_i\ft{(\tiK_\bul,y_\bul)}:=C_i\ft{(\tiK_\bul,y_\bul)}\setminus %
 C_{i-1}\ft{(\tiK_\bul,y_\bul)}=w^{\tiK_{i-1}}y_{i-1} %
 \vPhi_{\tiK_{i-1}}^\infty([\bs_{i-1}]), \]
 where $y_0=1$.

 Step2.
 For each $i=1,\dots,n$, define a total order
 $\preq_i\!\!\ft{(\tiK_\bul,y_\bul,\bs_\bul)}=\preq_i$ on
 $R_i\ft{(\tiK_\bul,y_\bul)}$ by setting
\[ w^{\tiK_{i-1}}y_{i-1}\phi_{\bs_{i-1}}(p) \,\preq_i\, %
 w^{\tiK_{i-1}}y_{i-1}\phi_{\bs_{i-1}}(q) \]
 for each $p\leq q$.

 Step 3.
 Define $\preq_\bul$ by extending $\preq_1,\dots,\preq_n$ to
 $\De_{\tiJ}\ft{(w,-)}=R_1\ft{(\tiK_\bul,y_\bul)} %
 \amalg\cdots\amalg R_n\ft{(\tiK_\bul,y_\bul)}$ as follows:
\[ R_i\ft{(\tiK_\bul,y_\bul)} %
 \,\preq_\bul\, R_{i'}\ft{(\tiK_\bul,y_\bul)} \]
 for each $i<i'$.
\end{defn}

\begin{thm} 
 For each non-empty subset ${\bd J}\subset{\SI}$ and
 each element $w\in{\sW}_{\tiJ}$, the assignment
 $(\scK_\bul,y_\bul,\bs_\bul)\mapsto\, %
 \preq_\bul\!\!\ft{(\tiK_\bul,y_\bul,\bs_\bul)}$
 defines a bijective mapping between
 $\bol{\cl T}\ft{(\scJ,w)}$ and
 $\bol{\cl{CO}}\De_{\tiJ}\ft{(w,-)}$,
 which maps $\bol{\cl T}_n\ft{(\scJ,w)}$ onto
 $\bol{\cl{CO}}_n\De_{\tiJ}\ft{(w,-)}$
 for each $n=1,\dots,\sharp{\bd J}$.
\end{thm}

\begin{pf}
 We first prove that
 $\preq_\bul\!\!\ft{(\tiK_\bul,y_\bul,\bs_\bul)}=\preq_\bul$
 is a $n$-row type convex order on $\De_{\tiJ}\ft{(w,-)}$ for
 each element $(\scK_\bul,y_\bul,\bs_\bul)$ of
 $\bol{\cl T}_n\ft{(\scJ,w)}$. By definition, we have
\[ C_i\ft{(\tiK_\bul,y_\bul)}=C_{i-1}\ft{(\tiK_\bul,y_\bul)} %
 \amalg R_i\ft{(\tiK_\bul,y_\bul)},\quad %
 R_i\ft{(\tiK_\bul,y_\bul)}=w^{\tiK_{i-1}}y_{i-1} %
 \vPhi_{\tiK_{i-1}}^\infty([\bs_{i-1}]) \]
 for $i=1,\dots,n$. By Proposition 3.5(2),
 we see that $\preq_1$ is a 1-row type convex order on
 $C_1\ft{(\tiK_\bul,y_\bul)}=R_1\ft{(\tiK_\bul,y_\bul)}$ and that
 there exists a unique 1-row type convex order $\preq_i'$ on
 $\vPhi_{\tiK_{i-1}}^\infty([\bs_{i-1}])$ satisfying
 $\phi_{\preq_i'}=\phi_{\bs_{i-1}}$ for each $i=2,\dots,n$.
 By applying Proposition 7.7(1) recursively, we see that
 the restriction of $\preq_\bul$ to $C_i\ft{(\tiK_\bul,y_\bul)}$
 is a convex order for each $i=1,\dots,n$.
 Thus $\preq_\bul$ is a $n$-row type convex order on
 $\De_{\tiJ}\ft{(w,-)}$, since
 $C_n\ft{(\tiK_\bul,y_\bul)}=\De_{\tiJ}\ft{(w,-)}$.

 The injectivity follows from Theorem 6.9 and Proposition 3.5(2).
 Let us check the surjectivity.
 Suppose that
 $\preq$ is a $n$-row type convex order on $\De_{\tiJ}\ft{(w,-)}$.
 Then there exists a unique
 $C_\bul\ft{(\preq)}\in\bol{\cl C}_n\De_{\tiJ}\ft{(w,-)}$
 with the properties (i) and (ii) in Proposition 7.5.
 By Theorem 6.9, there exists a unique
 $(\scK_\bul,y_\bul)\in\bol{\cl S}_n\ft{(\scJ,w)}$ such that
 $C_\bul\ft{(\preq)}=C_\bul\ft{(\tiK_\bul,y_\bul)}$.
 Moreover, by Proposition 6.10, there exists a unique
 $B_{i-1}\ft{(\tiK_\bul,y_\bul)}\in{\fr B}_{\tiK_{i-1}}^\infty$
 such that
\[ R_i\ft{(\preq)}=R_i\ft{(\tiK_\bul,y_\bul)}= %
 w^{\tiK_{i-1}}y_{i-1}B_{i-1}\ft{(\tiK_\bul,y_\bul)} \]
 for each $i=1,\dots,n$.
 Since the restriction of $\preq$ to $C_i\ft{(\preq)}$ is
 a convex order and the restriction $\preq_i$ of $\preq$ to
 $R_i\ft{(\preq)}$ is of 1-row type, there exists a unique 1-row type
 convex order $\preq_i'$ on $B_{i-1}\ft{(\tiK_\bul,y_\bul)}$ such that
 $\phi_{\ti{\preq_i}}=w^{\tiK_{i-1}}y_{i-1}\phi_{\ti{\preq_i'}}$
 by Proposition 7.7(2).
 Thus there exists a unique
 $\bs_{i-1}\in{\cl W}_{\tiK_{i-1}}^\infty$ such that
 $B_{i-1}\ft{(\tiK_\bul,y_\bul)}= %
 \vPhi_{\tiK_{i-1}}^\infty([\bs_{i-1}])$ and
 $\phi_{\ti{\preq_i}}=w^{\tiK_{i-1}}y_{i-1}\phi_{\bs_{i-1}}$
 by Proposition 3.5(2).
 Put
 $\bs_\bul:=(\bs_0,\dots,\bs_{n-1})\in{\cl W}_{\tiK_\bul}^\infty$.
 Then we have
 $(\scK_\bul,y_\bul,\bs_\bul)\in\bol{\cl T}_n\ft{(\scJ,w)}$ and
 $\preq_i=\preq_i\!\!\ft{(\tiK_\bul,y_\bul,\bs_\bul)}$
 for each $i=1,\dots,n$, and hence
 $\preq=\preq_\bul\!\!\ft{(\tiK_\bul,y_\bul,\bs_\bul)}$.
\end{pf}

\noindent
{\it Example\/.}
 In case $\De$ is the root system of affine Lie algebra of
 type $A_2^{(1)}$, we have ${\SI}=\{1,2\}$ and
 ${\sDe}_\pm=\pm\{\al_1,\al_1+\al_2,\al_2\}$.
 Let ${\bd K}_\bul$ be an element of $\bol{\cl C}_2{\SI}$
 such that ${\bd K}_1=\{1\}$, and $y_\bul$ an element of
 $W_{\tiK_\bul}$ such that $y_1=s_1s_{\de-\al_1}$.
 Further, let $\bs_0$ be an element of ${\cl W}_{\tiK_0}^\infty$
 such that
 $(\bs_0(1),\bs_0(2),\bs_0(3),\bs_0(4),\bs_0(5),\bs_0(6))= %
 (s_2,s_1,s_0,s_2,s_1,s_2)$ and
 $(\bs_0(3q+4),\bs_0(3q+5),\bs_0(3q+6))=(s_0,s_1,s_2)$
 for each $q\in{\bd N}$, and $\bs_1$ an element of
 ${\cl W}_{\tiK_1}^\infty$ such that
 $(\bs_1(2q-1),\bs_1(2q))=(s_1,s_{\de-\al_1})$ for each $q\in{\bd N}$.
 Put $\bs_\bul:=(\bs_0,\bs_1)\in{\cl W}_{\tiK_\bul}^\infty$.
 Then $(\scK_\bul,y_\bul,\bs_\bul)$ is an element of
 $\bol{\cl T}_2\ft{(\scri{\SI},w)}$, where $w=s_2s_1\in{\sW}$.
 The following total order $\preq$ on
 $\De\ft{(w,-)}$ is the convex order
 $\preq_\bul\!\!\ft{(\tiK_\bul,y_\bul,\bs_\bul)}$ of 2-row type:
\begin{align*}
 \al_2 \pr \al_1+\al_2 &\pr \de+\al_2 \pr \de+\al_1+\al_2 %
  \pr2\de+\al_2 \pr \de-\al_1 \\
 \pr 3\de+\al_2 &\pr 2\de-\al_1 \pr 4\de+\al_2 \pr 3\de-\al_1 %
  \pr 5\de+\al_2 \pr 4\de-\al_1 \pr\cdots\cdots \\
 \pr 2\de+\al_1+\al_2 &\pr 3\de+\al_1+\al_2 \pr 4\de+\al_1+\al_2 %
  \pr 5\de+\al_1+\al_2 \pr\cdots\cdots
\end{align*}
 with
\begin{align*}
 R_1\ft{(\preq)} &= %
 \{\,n\de-\al_1,\,m\de+\al_2,\,k\de+\al_1+\al_2\;|\; %
 n\in{\bd N},\,m\in{\bd Z}_{\geq0},\,k=0,1\,\}, \\
 R_2\ft{(\preq)} &= %
 \{\,k\de+\al_1+\al_2\;|\;k\in{\bd Z}_{\geq2}\,\}.
\end{align*}

\begin{cor} 
 Let us choose an element $w\in{\sW}$, an element
 $(\scJ_\bul,y_\bul,\br_\bul)\in\bol{\cl T}\ft{(\scri{\SI},w)}$,
 a total order $\preq_0$ on $\pim$, and an element
 $(\scK_\bul,z_\bul,\bs_\bul)\in %
 \bol{\cl T}\ft{(\scri{\SI},ww_\circ)}$, where
 $w_\circ\in{\sW}$ is the longest element.
 Define a total order $\preq$ on $\De_+$ by extending
 the convex order
 $\preq_\bul\!\!\ft{(\scJ_\bul,y_\bul,\br_\bul)}$ on
 $\De\ft{(w,-)}$, the total order $\preq_0$ on $\pim$,
 and the opposite convex order
 $\preq_\bul\!\!\ft{(\scK_\bul,z_\bul,\bs_\bul)}^{op}$ on
 $\De\ft{(w,+)}$ to
 $\De_+=\De\ft{(w,-)}\amalg\pim\amalg\De\ft{(w,+)}$
 as follows{\em:}
\[ \De\ft{(w,-)}\,\pr\,\pim\,\pr\,\De\ft{(w,+)}. \]
 Then $\preq$ is a convex order on $\De_+$.
 Moreover, each convex order on $\De_+$ is constructed by
 applying the above procedure suitably.
\end{cor}

\begin{pf}
 This follows immediately from Theorem 7.3 and Theorem 7.9.
\end{pf}